\documentclass[12pt]{article}

\usepackage{authblk}
\usepackage{pgfplots}
\usepackage{tikz}

\usetikzlibrary{automata}
\usetikzlibrary{positioning}
\usepackage{graphicx}

\usepackage{amsmath,amsthm,enumerate,graphics,amsfonts,latexsym,amsopn,verbatim,amscd,amssymb}
\usepackage[shortlabels]{enumitem}

\usepackage{txfonts,geometry,hyperref}

\newtheorem{conj}{Conjecture}[section]
\newtheorem{remark}{Remark}[section]

\newtheorem{theorem}{Theorem}[section]

\newtheorem{cor}[theorem]{Corollary}

\begin{document}
\title{Zagreb indices of commuting and non-commuting graphs of finite groups and Hansen-Vuki{\v{c}}evi{\'c} conjecture}
\author{Shrabani Das,    Arpita Sarkhel and Rajat Kanti Nath\footnote{corresponding author}}
\affil{Department of Mathematical Sciences, Tezpur University, Napaam-784028, Sonitpur, Assam, India.
	
Emails: shrabanidas904@gmail.com (S. Das);  arpitasarkhel17@gmail.com (A. Sarkhel);  rajatkantinath@yahoo.com (R. K. Nath)}
\date{ }
\maketitle


\begin{abstract}
In this paper we compute first and second Zagreb indices of  commuting and non-commuting graphs of  finite  groups and determine several classes of finite groups such that their commuting  and non-commuting graphs satisfy Hansen-Vuki{\v{c}}evi{\'c} conjecture.
\end{abstract}



\section{Introduction}

Let $\mathfrak{G}$ be the set of all graphs. A topological index  is a function $T :  \mathfrak{G} \to \mathbb{R}$ such that $T(\Gamma_1) = T(\Gamma_2)$ whenever the graphs $\Gamma_1$ and   $\Gamma_2$ are isomorphic.  By using different parameters of graphs many topological indices have been defined since 1947. Wiener index is the first topological index, introduced by Wiener \cite{Wiener-97}, and it is a distance based index. Among the degree based topological indices, the first two (known as Zagreb indices) were introduced by Gutman and Trinajsti{\'c} \cite{Gut-Trin-72} in 1972. Initially, topological indices were used to describe several chemical properties of molecules. In particular, Zagreb indices were used in examining the dependence of total $\pi$-electron energy on molecular structure. As noted in \cite{Z-index-30y-2003} , Zagreb indices are also used in studying molecular complexity, chirality,
ZE-isomerism and heterosystems etc. Later on, general mathematical properties of various topological indices  are also studied by many mathematicians. A survey on mathematical properties of Zagreb indices can be found in \cite{Gut-Das-2004}. Computing formulas for Zagreb indices of non-commuting graph $\mathcal{NC}(G)$ of a finite group $G$ were obtained in \cite{mirzargar2012some}. However, their formulas are not closed because of the presence of terms like $\sum_{x \in G \setminus Z(G)}C_G(x)$ and $\sum_{xy \in \ e(\mathcal{NC}(G))}|C_G(x)||C_G(y)|$, where $C_G(x) = \{g \in G : xg = yg\}$ (the centralizer of $x \in G$) and $e(\mathcal{NC}(G)$ is the set of edges of $\mathcal{NC}(G)$. Zagreb indices of commuting graphs  of groups are yet not  explored.

Let $\Gamma$ be a simple undirected graph with vertex set $v(\Gamma)$ and edge set $e(\Gamma)$. The first and second Zagreb indices of $\Gamma$, denoted by $M_{1}(\Gamma)$ and $M_{2}(\Gamma)$ respectively, are defined as 
\[
M_{1}(\Gamma) = \sum\limits_{v \in v(\Gamma)} \deg(v)^{2}  \text{ and }  M_{2}(\Gamma) = \sum\limits_{uv \in e(\Gamma)} \deg(u)\deg(v),
\]
where $ \deg(v) $ is the number of edges incident on $ v $
(called degree of $v$). Comparing first and second Zagreb indices, Hansen and Vuki{\v{c}}evi{\'c} \cite{hansen2007comparing} posed the following   conjecture in  2007.
\begin{conj}\label{Conj}
(Hansen-Vuki{\v{c}}evi{\'c} Conjecture) For any simple finite graph $\Gamma$, 
\begin{equation}\label{Conj-eq}
 \dfrac{M_{2}(\Gamma)}{\vert e(\Gamma) \vert} \geq \dfrac{M_{1}(\Gamma)}{\vert v(\Gamma) \vert} .
\end{equation}
\end{conj}
It was shown in \cite{hansen2007comparing} that the conjecture is not true  if $\Gamma = K_{1, 5} \sqcup K_3$. However, Hansen and Vuki{\v{c}}evi{\'c} \cite{hansen2007comparing}  showed that  Conjecture \ref{Conj}  holds for chemical graphs. In \cite{vukicevic2007comparing}, it was shown that the conjecture holds for trees with equality in \eqref{Conj-eq} when $\Gamma$ is a star graph. In \cite{liu2008conjecture}, it was shown that the conjecture holds for connected unicyclic graphs with equality when the graph is a cycle. The case when equality holds in \eqref{Conj-eq} is studied extensively in \cite{vukicevic2011some}. A survey on comparing Zagreb indices can be found in \cite{Liu-You-2011}. Interestingly, it is not known whether Conjecture \ref{Conj} holds for commuting and non-commuting graphs of finite groups. In this paper, we compute first and second Zagreb indices of commuting and non-commuting graphs of several families of finite non-abelian groups and check the validity of Hansen-Vuki{\v{c}}evi{\'c} Conjecture. It is worth mentioning that  Zagreb indices of commuting conjugacy class graph and its complement are computed and verified  Conjecture \ref{Conj} in \cite{Das-Nath-2023} for the classes of finite groups considered in \cite{Salah-2020, SA-2020, SA-CA-2020}.
 
The commuting graph $\mathcal{C}(G)$ of a finite non-abelian group $G$ is a graph defined on the elements of $G \setminus Z(G)$ and two elements $x$ and $y$ are adjacent if and only if $ xy=yx$. The complement of $\mathcal{C}(G)$ (also denoted by $\overline{\mathcal{C}(G)}$ ) is nothing but $\mathcal{NC}(G)$. 
 The commuting graph was first studied by Brauer and Fowler \cite{MR74414}, in the year 1955.  For the structures of commuting graphs of various classes of finite non-abelian groups we refer \cite{DBN-KJM-2020,B1,B2,DN-KJM-2018,  DN-IJPAM-2021,FSN-ADM-2021,NFDS-SB-2021,SND-PJM-2022}, where various spectra and energies of $C(G)$  were computed.  
 
\section{Zagreb indices of commuting and non-commuting graphs}
In this section, we consider several classes of well-known finite groups and compute  Zagreb indices of their commuting and non-commuting graphs. 
The following results are useful in the computations.
\begin{theorem}\label{thm1}
Let $\Gamma$ be the disjoint union of the graphs  $ \Gamma_{1}, \Gamma_{2}, \dots, \Gamma_{n}$. If $\Gamma_{i} = l_{i}K_{m_{i}}$ for $ i = 1, 2, \dots, k$, where $K_{m_{i}}$'s are complete graphs on $ m_{i} $ vertices and $l_iK_{m_i}$ is the disjoint union of $l_i$ copies of $K_{m_i}$, then
\[ M_{1}(\Gamma) = \sum_{i=1}^{k}l_{i}m_{i}(m_{i}-1)^{2} ~~~ \text{and} ~~~ M_{2}(\Gamma) = \sum_{i=1}^{k}l_{i}\dfrac{m_{i}(m_{i}-1)^{3}}{2}.\]
\end{theorem}
\begin{proof}
    By definitions of $M_{1}(\Gamma)$ and $ M_{2}(\Gamma)$ we have
\begin{equation}\label{Comp-eq-1}
M_{1}(\Gamma) = \sum_{i = 1}^{k} M_{1}(\Gamma_i) \quad \text{and} \quad M_{2}(\Gamma) = \sum_{i = 1}^{k} M_{2}(\Gamma_i).
\end{equation}
    If $\Gamma_{i} = l_{i}K_{m_{i}}$ for $ i = 1, 2, \dots, k$ then 
\begin{equation}\label{Comp-eq-2}
M_{1}(\Gamma_i) = l_iM_{1}(K_{m_i}) \quad \text{and} \quad M_{2}(\Gamma) = l_i M_{2}(K_{m_i}).
\end{equation}
    Hence, the result follows from \eqref{Comp-eq-1} and \eqref{Comp-eq-2} noting that
\[ 
M_{1}(K_{m_i}) = m_i(m_i-1)^{2} ~~~ \text{and} ~~~ M_{2}(K_{m_i}) = \frac{m_i(m_i-1)^{3}}{2}.
\]
\end{proof}
\begin{theorem}\label{using commuting finding complement}
    (\cite{Das-Kex-14}, Page 575 and \cite{Das-Gut-04}, Lemma 3) For any graph  $\Gamma$ and its complement $\overline{\Gamma}$,
\[
M_{1}(\overline{\Gamma})=|v(\Gamma)|(|v(\Gamma)|-1)^{2}-4|e(\Gamma)|(|v(\Gamma)|-1)+ M_{1}(\Gamma) \quad  \text{ and }
\]
\[ 
M_{2}(\overline{\Gamma})= \dfrac{|v(\Gamma)|(|v(\Gamma)|-1)^{3}}{2}+2|e(\Gamma)|^{2}-3|e(\Gamma)|(|v(\Gamma)|-1)^{2}+\left(|v(\Gamma)|-\dfrac{3}{2}\right)M_{1}(\Gamma)-M_{2}(\Gamma). 
\]
\end{theorem}

We first consider $\mathcal{C}(G)$ and $\mathcal{NC}(G)$ for the groups $G = D_{2m}, Q_{4n}, QD_{2^n}, SD_{8n}$ and $V_{8n}$.
\begin{theorem}\label{thm-D2m}
    If $G=D_{2m} = \langle f, g :  f^m = g^2 = 1, gfg^{-1} = f^{-1} \rangle$ $(m \geq 3)$, then 
    $$
    M_{1}(\mathcal{C}(G)) = \begin{cases}
        (m-1)(m-2)^{2}, & \text{when $m$ is odd} \\
        (m-2)(m-3)^{2}+m, & \text{when $m$ is even,}
    \end{cases} 
    $$
    $$
    M_{2}(\mathcal{C}(G))= \begin{cases}
        \dfrac{(m-1)(m-2)^{3}}{2}, & \text{when $m$ is odd} \\
        \dfrac{(m-2)(m-3)^{3}}{2} + \dfrac{m}{2}, & \text{when $m$ is even,} 
    \end{cases}
    $$
    $$ 
    M_{1}(\mathcal{NC}(G)) = \begin{cases}
        m(m-1)(5m-4), & \text{when m is odd} \\
        5m^{3}-18m^{2}+16m, & \text{when m is even}
    \end{cases}
    $$  
and
    $$ 
    M_{2}(\mathcal{NC}(G))= \begin{cases}
        m(m-1)(4m^{2}-6m+2), ~ \text{when m is odd} \\
       4m^{4}-20m^{3}+32m^{2}-16m, ~ \text{when m is even} 
    \end{cases} 
    $$
    Further, $\frac{M_{2}(\Gamma(G))}{|e(\Gamma(G))|} \geq \frac{M_{1}(\Gamma(G))}{|v(\Gamma(G))|}$, where $\Gamma(G) = \mathcal{C}(G)$ or $\mathcal{NC}(G)$, with equality when $m=4$.   
\end{theorem}
\begin{proof}
    \textbf{Case 1.} $ m $ is odd.

    It is well-known that $\mathcal{C}(D_{2m}) = K_{m-1} \sqcup mK_{1} $.
    As such, $ |v( \mathcal{C}(D_{2m}) )| = 2m-1$ and $|e(\mathcal{C}(D_{2m}) )| = \binom{m-1}{2}= \dfrac{(m-1)(m-2)}{2} $. Therefore, using Theorem \ref{thm1}, we get
\begin{align*}
    M_{1}(\mathcal{C}(D_{2m})) &=  (m-1)(m-1-1)^{2} + m(1-1)^{2}
    = (m-1)(m-2)^{2} \quad \text{ and }
\end{align*}  
\begin{align*}
    M_{2}(\mathcal{C}(D_{2m})) = \dfrac{(m-1)(m-1-1)^{3}}{2} + m \cdot \dfrac{1(1-1)^{3}}{2} 
     =\dfrac{(m-1)(m-2)^{3}}{2}.
\end{align*} 
We have
\[ 
\dfrac{M_{1}(\mathcal{C}(D_{2m}))}{|v( \mathcal{C}(D_{2m}) )|} = \frac{(m-1)(m-2)^{2}}{2m-1} \quad  \text{ and } \quad \dfrac{M_{2}(\mathcal{C}(D_{2m}))}{|e(\mathcal{C}(D_{2m}))|} = (m-2)^{2}.
\]
    Also, for $m \geq 3$ we have $m-1  < 2m - 1$ and so  $(m-2)^{2}  >  \frac{(m-1)(m-2)^{2}}{2m-1}$. Therefore,
\[
\dfrac{M_{2}(\mathcal{C}(D_{2m}))}{|e(\mathcal{C}(D_{2m}))|}  > \dfrac{M_{1}(\mathcal{C}(D_{2m}))}{|v(\mathcal{C}(D_{2m}))|}.
\]
    Using Theorem \ref{using commuting finding complement} we have 
 \begin{align*}
      M_{1}(\mathcal{NC}(D_{2m}))&=(2m-1)(2m-2)^{2}-4(2m-2)\dfrac{(m-1)(m-2)}{2}+(m-1)(m-2)^{2} \\
      &= (m-1)[(8m-4)(m-1)-4(m-1)(m-2)+(m-2^{2})]\\
      &= m(m-1)(5m-4)
\end{align*} 
 and 
\begin{align*}
      M_{2}(\mathcal{NC}(D_{2m}))&=\dfrac{(2m-1)(2m-2)^{3}}{2}+2\dfrac{(m-1)^{2}(m-2)^{2}}{4}-3 \dfrac{(m-1)(m-2)}{2}(2m-2)^{2} \\
     &~~~~~~~~~~~~~~~~~~~~~~~~~~~~~~~~~~~~~~~~~~~~~~+(2m-1-\dfrac{3}{2})(m-1)(m-2)^{2}-\dfrac{(m-1)(m-2)^{3}}{2} \\
      &= \dfrac{m-1}{2}(8m^{3}-12m^{2}+4m) \\
      &= m(m-1)(4m^{2}-6m+2).
\end{align*}
    Also, $ |v(\mathcal{NC}(D_{2m}))|=2m-1 $ and $ |e(\mathcal{NC}(D_{2m}))| = \binom{2m-1}{2}-|e(\mathcal{C}(D_{2m}))|=\frac{3m(m-1)}{2} $. We have 
 \[
 \dfrac{M_{1}(\mathcal{NC}(D_{2m}))}{|v( \mathcal{NC}(D_{2m}) )|} = \dfrac{m(m-1)(5m-4)}{2m-1} \text{ and } \dfrac{M_{2}(\mathcal{NC}(D_{2m}))}{|e( \mathcal{NC}(D_{2m}) )|} =  \dfrac{m(m-1)(8m^{2}-12m+4)}{3m(m-1)}. 
 \]
 As such
 \[
 \dfrac{M_{2}(\mathcal{NC}(D_{2m}))}{|e( \mathcal{NC}(D_{2m}) )|} - \dfrac{M_{1}(\mathcal{NC}(D_{2m}))}{|v( \mathcal{NC}(D_{2m}) )|} = \dfrac{m^{2}(m-5)+4(2m-1)}{3m(m-1)(2m-1)}:= \dfrac{f(m)}{g(m)}. 
 \]
 Since  $f(m), g(m) > 0$   for all $m \geq 3$ we have $\frac{f(m)}{g(m)} > 0$.
 
 \noindent \textbf{Case 2.} $ m $ is even.

    It is well-known that  $\mathcal{C}(D_{2m}) = K_{m-2} \sqcup \frac{m}{2}K_{2} $. As such, $|v( \mathcal{C}(D_{2m}) )|= 2m-2 $ and ~ $ |e(\mathcal{C}(D_{2m}))| = \binom{m-2}{2} + \frac{m}{2}= \frac{(m-2)(m-3) + m}{2}$. Therefore, using Theorem \ref{thm1}, we get 
\begin{align*}
    M_{1}(\mathcal{C}(D_{2m})) = (m-2)(m-2-1)^{2} + \dfrac{m}{2} \cdot 2(2-1)^{2} 
    = (m-2)(m-3)^{2} + m  \quad \text{ and }
\end{align*} 
 \begin{align*}
     M_{2}(\mathcal{C}(D_{2m})) = \dfrac{(m-2)(m-2-1)^{3}}{2} + \dfrac{m}{2} \cdot \dfrac{2(2-1)^{3}}{2} 
     = \dfrac{(m-2)(m-3)^{3} + m}{2}.
\end{align*} 
We have
\[
\dfrac{M_{2}(\mathcal{C}(D_{2m}))}{|e(\mathcal{C}(D_{2m}))|} = \dfrac{(m-2)(m-3)^{3}+m}{(m-2)(m-3)+m} \quad  \text{ and } \quad  
\dfrac{M_{1}(\mathcal{C}(D_{2m}))}{|v( \mathcal{C}(D_{2m}) )|} = \dfrac{(m-2)(m-3)^{2}+m}{2m-2}.
\]
For $m=4$ we have  
\[
\dfrac{M_{2}(\mathcal{C}(D_{2m}))}{|e(\mathcal{C}(D_{2m}))|} = 1 =  
\dfrac{M_{1}(\mathcal{C}(D_{2m}))}{|v( \mathcal{C}(D_{2m}))|}.
\]
For $m \geq 6$ we have
\[
(m-3)^{3}+1 - (m-3) - (m-3)^{2} = (m-3)((m-3)(m-4) - 1) + 1 > 0.
\]
Therefore,
\[
(m-3)^{3}+1 > (m-3) + (m-3)^{2}.
\]
Multiplying both sides by $m(m-2)$ we get
\[
(m-2)(m-3)^{3}m + m(m-2) >  m(m-2)(m-3) + m(m-2)(m-3)^{2}.
\]
Adding $(m-2)^2(m-3)^3 + m^2$ we get
\begin{align*}
(m-2)(m-3)^{3}(2m-2) &+ m(2m-2)\\
 & > (m-2)^{2}(m-3)^{3} + m(m-2)(m-3) + m(m-2)(m-3)^{2} + m^{2} \\
& = ((m-2)(m-3)+m)((m-2)(m-3)^{2}+m).
\end{align*}
Therefore,
\[
\dfrac{M_{2}(\mathcal{C}(D_{2m}))}{|e(\mathcal{C}(D_{2m}))|} = \dfrac{(m-2)(m-3)^{3}+m}{(m-2)(m-3)+m} > \dfrac{(m-2)(m-3)^{2}+m}{2m-2} = \dfrac{M_{1}(\mathcal{C}(D_{2m}))}{|v( \mathcal{C}(D_{2m}))|}.
\]
    Using Theorem \ref{using commuting finding complement} we have 
\begin{align*}
      M_{1}(\mathcal{NC}(D_{2m}))&=(2m-2)(2m-3)^{2}-4(2m-3)\dfrac{(m-2)(m-3)+m}{2}+(m-2)(m-3)^{2}+m \\
      &=  5m^{3}-18m^{2}+16m
\end{align*} 
and 
\begin{align*}
      M_{2}(\mathcal{NC}(D_{2m}))&=\dfrac{(2m-2)(2m-3)^{3}}{2}+2\dfrac{(m-2)^{2}(m-3)^{2}+2m(m-2)(m-3)+m^{2}}{4} \\
      &~~~~~~~~~~-3 \dfrac{(m-2)(m-3)+m}{2}(2m-3)^{2}+(2m-2-\dfrac{3}{2})((m-2)(m-3)^{2}+m) \\
     &~~~~~~~~~~~~~~~~~~~~~~~~~~~~~-\dfrac{(m-2)(m-3)^{3}+m}{2} \\
      &= \dfrac{1}{2}(8m^{4}-40m^{3}+64m^{2}-32m) \\
      &=4m^{4}-20m^{3}+32m^{2}-16m .
 \end{align*}
    Also, $ |v(\mathcal{NC}(D_{2m}))|=2m-2 $ and $ |e( \mathcal{NC}(D_{2m}) )| = \binom{2m-2}{2}-|e(\mathcal{NC}(D_{2m}))|=\frac{3m(m-2)}{2} $. We have 
 \[
 \dfrac{M_{1}(\mathcal{NC}(D_{2m}))}{|v( \mathcal{NC}(D_{2m}))|} = \dfrac{5m^{3}-18m^{2}+16m}{2m-2} 
 \]
 and
 \[
 \dfrac{M_{2}(\mathcal{NC}(D_{2m}))}{|e( \mathcal{NC}(D_{2m}))|} =  \dfrac{8m^{4}-40m^{3}+64m^{2}-32m}{3m(m-2)}. 
 \]
    As such 
 \[
 \dfrac{M_{2}(\mathcal{NC}(D_{2m}))}{|e( \mathcal{NC}(D_{2m}) )|} - \dfrac{M_{1}(\mathcal{NC}(D_{2m}))}{|v( \mathcal{NC}(D_{2m}) )|} = \dfrac{m^{3}(m^{2}-12m+28)+m^{2}(24m-96)+64m}{3m(m-2)(2m-2)}:= \frac{f(m)}{g(m)}. 
 \]
 We have $ g(m) > 0$   for all $m \geq 4$, $f(4)=0$, $f(6)=384$ and $f(8)=4608$. For $m \geq 10$ we have $m^{2}-12m+28 > 0$, $24m-96>0$ and so $f(m)>0$. Hence $\frac{f(m)}{g(m)} \geq 0$ with equality when $m=4$.
\end{proof}
\begin{cor}\label{cor-Q4n}
    If $ G=Q_{4n} = \langle f, g : f^{2n} = 1, g^2 = f^n, gfg^{-1} = f^{-1} \rangle$ $(n \geq 2)$, then 
\[
 M_{1}(\mathcal{C}(G)) =(2n-2)(2n-3)^{2}+2n, \quad M_{2}(\mathcal{C}(G))= (n-1)(2n-3)^{3}+n,
 \]
 \[
 M_{1}(\mathcal{NC}(G))=40n^{3}-72n^{2}+32n \quad \text{ and } \quad M_{2}(\mathcal{NC}(G))=64n^{4}-160n^{3}+128n^{2}-32n.
\]
    Further, $\frac{M_{2}(\Gamma(G))}{|e(\Gamma(G))|} \geq \frac{M_{1}(\Gamma(G))}{|v(\Gamma(G))|}$, where $\Gamma(G) = \mathcal{C}(G)$ or $\mathcal{NC}(G)$, with equality if and only if $n=2$.
\end{cor}
\begin{proof}
    It is well-known that  $\mathcal{C}(Q_{4n})= K_{2n-2} \sqcup nK_{2} \cong \mathcal{C}({D_{2 \times 2n}})$. Therefore, putting $m=2n$ in  Theorem \ref{thm-D2m}, we get the required result.
\end{proof}
\begin{cor}\label{cor-QD2n}
    If $ G=QD_{2^{n}} = \langle f, g : f^{2^n} = g^2 = 1, gfg^{-1} = f^{-1} \rangle$ ($n\geq 3$), then 
\[
M_{1}(\mathcal{C}(G)) = (2^{n-1}-2)(2^{n-1}-3)^{2} + 2^{n-1}, \quad M_{2}(\mathcal{C}(G))= (2^{n-2}-1)(2^{n-1}-3)^{3} + 2^{n-2},
\]
\[
M_{1}(\mathcal{NC}(G))=5 \cdot 2^{3n-3}-18 \cdot 2^{2n-2} +16 \cdot 2^{n-1} \quad \text{ and } 
\]
\[
M_{2}(\mathcal{NC}(G))=4 \cdot 2^{4n-4}-20 \cdot 2^{3n-3}+32 \cdot 2^{2n-2}- 16 \cdot 2^{n-1}.
\]
Further, $\frac{M_{2}(\Gamma(G))}{|e(\Gamma(G))|} \geq \frac{M_{1}(\Gamma(G))}{|v(\Gamma(G))|}$, where $\Gamma(G) = \mathcal{C}(G)$ or $\mathcal{NC}(G)$, with equality if and only if $n=3$.
\end{cor}
\begin{proof}
    It is well-known that  $\mathcal{C}(QD_{2^{n}})= K_{2^{n-1}-2} \sqcup 2^{n-2}K_{2} \cong \mathcal{C}({D_{2 \times 2^{n-1}}})$. Therefore, putting $m=2^{n-1}$ in  Theorem \ref{thm-D2m}, we get the required result.
\end{proof}
\begin{theorem}\label{thm-V8n}
    If $G=V_{8n} = \langle f, g: f^{2n}=g^{4}=1, gf=g^{-1}f^{-1}, g^{-1}f=f^{-1}g\rangle$, then 
	$$
	M_{1}(\mathcal{C}(G))=\begin{cases}
		(4n-4)(4n-5)^{2}+36n, &  \text{when $n$ is even} \\
		(4n-2)(4n-3)^{2}+4n, & \text{when $n$ is odd,} 
	\end{cases}
	$$  
	$$
	M_{2}(\mathcal{C}(G))=  \begin{cases}
		(2n-2)(4n-5)^{3}+54n, & \text{when $n$ is even} \\
		(2n-1)(4n-3)^{3}+2n, & \text{when $n$ is odd,}
	\end{cases}
        $$
        $$
        M_{1}(\mathcal{NC}(G))= \begin{cases}
            8n(40n^{2}+8n-93), & \text{when n is even} \\
            16n(20n^{2}-18n+4), & \text{when n is odd}
        \end{cases}
        $$
and
        $$
        M_{2}(\mathcal{NC}(G))=\begin{cases}
           2n(512n^{3}-1180n^{2}+1024n-229), ~ \text{when n is even} \\
           64n(16n^{3}-20n^{2}+8n-1), ~ \text{when n is odd.}
        \end{cases}
        $$
	Further,  	$\frac{M_{2}(\Gamma(G))}{|e(\Gamma(G))|} \geq \frac{M_{1}(\Gamma(G))}{|v(\Gamma(G))|}$, where $\Gamma(G) = \mathcal{C}(G)$ or $\mathcal{NC}(G)$, with equality when $n=1, 2$.
\end{theorem}
\begin{proof}
	\textbf{Case 1.}  $n$ is even. 
	
    It is well-known that $\mathcal{C}(G)=K_{4n-4} \sqcup nK_{4}.$ As such, $|v(\mathcal{C}(G))| = 8n-4$ and $|e(\mathcal{C}(G))|=\binom{4n-4}{2}+n\cdot\binom{4}{2}=(2n-2)(4n-5)+6n.$  Therefore, using Theorem \ref{thm1}, we get 
\begin{align*}
    M_{1}(\mathcal{C}(G))=(4n-4)(4n-4-1)^{2}+n \cdot4(4-1)^{2} = (4n-4)(4n-5)^{2}+36n \quad \text{ and }
\end{align*}
\begin{align*}
    M_{2}(\mathcal{C}(G))= \dfrac{(4n-4)(4n-4-1)^{3}}{2}+n \cdot \dfrac{4(4-1)^{3}}{2}= (2n-2)(4n-5)^{3}+54n.
\end{align*}
    We have
\[
\dfrac{M_{1}(\mathcal{C}(G))}{|v(\mathcal{C}(G))|}= \dfrac{(4n-4)(4n-5)^{2}+36n}{8n-4} \quad
\text{  and } \quad
\dfrac{ M_{2}(\mathcal{C}(G))}{|e(\mathcal{C}(G))|}= \dfrac{(2n-2)(4n-5)^{3}+54n}{(2n-2)(4n-5)+6n}. 
\]
    Therefore,
\begin{equation}\label{new-eq-V8n}
\dfrac{ M_{2}(\mathcal{C}(G))}{|e(\mathcal{C}(G))|} - \dfrac{M_{1}(\mathcal{C}(G))}{|v(\mathcal{C}(G))|} = \dfrac{32n^{4}(2n-11) + 32n^{2}(21n - 16) + 128n}{8n^{2}(n-2)+16n-5} := \frac{f(n)}{g(n)}.
\end{equation} 
We have $g(n) > 0$ for all $n \geq 2$, $f(2) = 0$ and $f(4)=\frac{10752}{315}$.   For $n \geq 6$ we have $2n-11 > 0$, $21n - 16 > 0$  and so $f(n) > 0$. Hence, $\frac{f(n)}{g(n)} \geq 0$ with equality when $n = 2$.

    Using Theorem \ref{using commuting finding complement} we have
\begin{align*}
    M_{1}(\mathcal{NC}(G))&=(8n-4)(8n-5)^{2}-4(8n-5)((2n-2)(4n-5)+6n) \\ 
    &~~~~~~~~~~~~~~~~~~~~~~~~~~~~~~~~~~~~~~~~~~~~~~~~~~~~+(4n-4)(4n-5)^{2}+36n \\
    &= 320n^{3}-576n^{2}+256n \\
    &= 8n(40n^{2}-72n+32)
\end{align*} 
    and 
\begin{align*}
    M_{2}(\mathcal{NC}(G))&=\dfrac{(8n-4)(8n-5)^{3}}{2}+2((2n-2)^{2}(4n-5)^{2}+12n(2n-2)(4n-5)+36n^{2}) \\
    &~~~~~~-3((2n-2)(4n-5)+6n)(8n-5)^{2} +(8n-4-\dfrac{3}{2})((4n-4)(4n-5)^{2}+36n)\\ 
    &~~~~~~~~~~~~~~-(2n-2)(4n-5)^{3}-54n \\
    &= 1024n^{4}-2560n^{3}+2048n^{2}-512n \\
    &= 2n(512n^{3}-1280n^{2}+1024n-256).
\end{align*}
    Also, $ |v(\mathcal{NC}(G) )|=8n-4 $ and $ |e( \mathcal{NC}(G) )| = \binom{8n-4}{2}-|e(\mathcal{C}(G))|=24n(n-1) $. We have
\[
\dfrac{M_{1}(\mathcal{NC}(G))}{|v(\mathcal{NC}(G))|} = \dfrac{8n(40n^{2}-72n+32))}{8n-4} \text{ and } \dfrac{M_{2}(\mathcal{NC}(G))}{|e( \mathcal{NC}(G) )|} =  \dfrac{2n(512n^{3}-1280n^{2}+1024n-256)}{24n(n-1)}.
\]
    As such 
\begin{equation}\label{new-V-8n-2}
\dfrac{M_{2}(\mathcal{NC}(G))}{|e( \mathcal{NC}(G))|} - \dfrac{M_{1}(\mathcal{NC}(G))}{|v( \mathcal{NC}(G))|} = \dfrac{64n^{3}(n-6)+ 64n(13n - 12) + 256}{24(n-1)(2n-1)}:= \dfrac{f(n)}{g(n)}.
\end{equation}
for $M_{1}(\mathcal{C}(G)), M_{2}(\mathcal{C}(G)), M_{1}(\mathcal{NC}(G))$ and $M_{2}(\mathcal{NC}(G))$. Further,
    We have $g(n) > 0$ for all $n \geq 2$, $f(2)=0$ and $ f(4)=2304 $. For $ n \geq 6$ we have $f(n) > 0 $. Therefore, $\frac{f(n)}{g(n)} \geq 0$ with equality when $n = 2$.
	
\noindent \textbf{Case 2.}  $n$ is odd.
	
    It is well-known that $\mathcal{C}(G)=K_{4n-2} \sqcup 2nK_{2} \cong \mathcal{C}({D_{2 \times 4n}})$. Therefore, putting $m=4n$ in  Theorem \ref{thm-D2m}, we get the required expressions for $M_{1}(\mathcal{C}(G)), M_{2}(\mathcal{C}(G)), M_{1}(\mathcal{NC}(G))$ and $M_{2}(\mathcal{NC}(G))$. Further, 	$\frac{M_{2}(\Gamma(G))}{|e(\Gamma(G))|} \geq \frac{M_{1}(\Gamma(G))}{|v(\Gamma(G))|}$, where $\Gamma(G) = \mathcal{C}(G)$ or $\mathcal{NC}(G)$, with equality when $n=1$.
\end{proof}

\begin{cor}
    If $G=SD_{8n} = \langle f, g : f^{4n} = g^2 = 1, gfg = f^{2n-1} \rangle$ $(n \geq 2) $, then
$$
M_{1}(\mathcal{C}(G))=\begin{cases}
    (4n-4)(4n-5)^{2}+36n , & \text{when n is odd} \\
    (4n-2)(4n-3)^{2}+4n, & \text{when n is even,} 
    \end{cases}
$$ 
$$
M_{2}(\mathcal{C}(G))=\begin{cases}
    (2n-2)(4n-5)^{3}+54n, & \text{when n is odd} \\
    (2n-1)(4n-3)^{3}+2n, & \text{when n is even,}
    \end{cases}
$$
$$
M_{1}(\mathcal{NC}(G))= \begin{cases}
    8n(40n^{2}+8n-93), & \text{when n is odd} \\
    16n(20n^{2}-18n+4), & \text{when n is even}
   \end{cases}
$$
$$ \text{and } \quad M_{2}(\mathcal{NC}(G))=\begin{cases}
    2n(512n^{3}-1180n^{2}+1024n-229), & \text{when n is odd} \\
   64n(16n^{3}-20n^{2}+8n-1), & \text{when n is even.}
   \end{cases}
$$
    Further, 	$\frac{M_{2}(\Gamma(G))}{|e(\Gamma(G))|} \geq \frac{M_{1}(\Gamma(G))}{|v(\Gamma(G))|}$, where $\Gamma(G) = \mathcal{C}(G)$ or $\mathcal{NC}(G)$, with equality when $n=2$. 
\end{cor}
\begin{proof}
\textbf{Case 1.} $n$ is odd.

   It is well-known that $\mathcal{C}(SD_{8n})=K_{4n-4} \sqcup nK_{4} $. Therefore, proceeding as in the proof of  Theorem \ref{thm-V8n} (Case 1)  we get the required expressions for $M_{1}(\mathcal{C}(G)), M_{2}(\mathcal{C}(G)), M_{1}(\mathcal{NC}(G))$, $M_{2}(\mathcal{NC}(G))$ and equations \eqref{new-eq-V8n} and \ref{new-V-8n-2}. Since $n \geq 3$ we have $\frac{M_{2}(\Gamma(G))}{|e(\Gamma(G))|} > \frac{M_{1}(\Gamma(G))}{|v(\Gamma(G))|}$.

\noindent \textbf{Case 2.} $n$ is even. 

  It is well-known that $\mathcal{C}(SD_{8n})=K_{4n-2} \sqcup 2nK_{2} \cong \mathcal{C}({D_{2 \times 4n}})$. Therefore, putting $m=4n$ in Theorem \ref{thm-D2m}, we get the required expressions for $M_{1}(\mathcal{C}(G)), M_{2}(\mathcal{C}(G)), M_{1}(\mathcal{NC}(G))$ and $M_{2}(\mathcal{NC}(G))$. Further,  
  $\frac{M_{2}(\Gamma(G))}{|e(\Gamma(G))|} \geq \frac{M_{1}(\Gamma(G))}{|v(\Gamma(G))|}$ with equality when $n=2$.
\end{proof}

Note that $\frac{G}{Z(G)}$ is isomorphic to some dihedral group if $G$ is itself a dihedral group or $G = Q_{4n}, QD_{2^n}$
and $SD_{8n}$ (when $n$ is even). This motivates us in obtaining the following result.

\begin{theorem}\label{thm-G/Z(G)=D_2m}
	Let $G$ be a finite group such that $\frac{G}{Z(G)} \cong D_{2m}, m \geq 3.$ Then  
	\[
	M_{1}(\mathcal{C}(G))= n(m-1)(mn-n-1)^{2} + mn(n-1)^{2}, M_{2}(\mathcal{C}(G))= \dfrac{(mn-n)(mn-n-1)^{3}+mn(n-1)^{3}}{2},
	\] 
	\[
	M_{1}(\mathcal{NC}(G))=n^{3}(5m^{3}-9m^{2}+4m) \quad \text{ and } \quad M_{2}(\mathcal{NC}(G))=n^{4}(4m^{4}-10m^{3}+8m^{2}-2m),
	\]
	where $n=|Z(G)|$. Further, 	$\frac{M_{2}(\Gamma(G))}{|e(\Gamma(G))|} > \frac{M_{1}(\Gamma(G))}{|v(\Gamma(G))|}$, where $\Gamma(G) = \mathcal{C}(G)$ or $\mathcal{NC}(G)$.
\end{theorem}
\begin{proof}
	It is well-known that $\mathcal{C}(G)=K_{(m-1)n} \sqcup mK_{n}$, where $n=|Z(G)|$. As such, $ |v(\mathcal{C}(G))|=  (2m-1)n$ and $|e(\mathcal{C}(G))| = \binom{mn-n}{2} + m \cdot \binom{n}{2} = \frac{(mn-n)(mn-n-1)+mn(n-1)}{2}.$ Therefore, using Theorem \ref{thm1}, we get
	$$
	M_{1}(\mathcal{C}(G))= n(m-1)(mn-n-1)^{2} + mn(n-1)^{2}  \quad \text{ and}
	$$ 
	\begin{align*}
		M_{2}(\mathcal{C}(G))&= \dfrac{(mn-n)(mn-n-1)^{3}}{2} + m \cdot \dfrac{n(n-1)^{3}}{2} \\
		&= \dfrac{(mn-n)(mn-n-1)^{3}+mn(n-1)^{3}}{2}.
	\end{align*}
	Also, 
	$$ 
	\dfrac{M_{1}(\mathcal{C}(G))}{|v(\mathcal{C}(G))|} = \dfrac{(m-1)(mn-n-1)^{2} + m(n-1)^{2}}{2m-1} 
	$$ 
	and 
	\[
	\dfrac{M_{2}(\mathcal{C}(G))}{|e(\mathcal{C}(G))|} = \dfrac{(m-1)(mn-n-1)^{3}+m(n-1)^{3}}{(m-1)(mn-n-1)+m(n-1)}.
	\]
	We have $(mn-2)^{2}-4(mn-n-1)(n-1) = mn^{2}(m-4)+4n^{2} > 0$. Therefore,
	\[
	(mn-2)^{2} - 3(mn-n-1)(n-1) > (mn-n-1)(n-1).
	\]
	Multiplying both sides by $(mn -2)$ we get
	\[
	(mn-2)^{3} - 3(mn-n-1)(n-1)(mn -2) > (mn-n-1)(n-1)(mn -2).
	\]
	We have  $(mn-2)^{3} - 3(mn-n-1)(n-1)(mn -2) = (mn-n-1)^3 + (n-1)^3$
	and so
	\[
	(mn-n-1)^3 + (n-1)^3 > (mn-n-1)(n-1)(mn -2).
	\]
	Multiplying both sides by $m(m-1)$ we get
	\[
	f(m, n) := m(m-1)(mn-n-1)^3 + m(m-1)(n-1)^3 > m(m-1)(mn-n-1)(n-1)(mn -2).
	\]
	Again,
	\[
	f(m, n) = (m-1)(2m - 1)(mn-n-1)^3 - (m - 1)^2(mn-n-1)^3 
	+ m(2m - 1)(n - 1)^3 - m^2(n-1)^2
	\]
	and so 
	\begin{align*}
		(m-1)&(2m - 1)(mn-n-1)^3 + m(2m - 1)(n - 1)^3\\
		& > (m - 1)^2(mn-n-1)^3 +  m^2(n-1)^2 + m(m-1)(mn-n-1)(n-1)(mn -2)\\
		& = ((m-1)(mn-n-1)+m(n-1))((m-1)(mn-n-1)^{2}+m(n-1)^{2}).
	\end{align*}
	Therefore,
	\[
	\dfrac{(m-1)(mn-n-1)^{3}+m(n-1)^{3}}{(m-1)(mn-n-1)+m(n-1)} > \dfrac{(m-1)(mn-n-1)^{2} + m(n-1)^{2}}{2m-1}
	\]
and so  $\frac{M_{2}(\mathcal{C}(G))}{|e(\mathcal{C}(G))|} > \frac{M_{1}(\mathcal{C}(G))}{|v(\mathcal{C}(G))|}$.	
	
	Using Theorem \ref{using commuting finding complement} we have 
	\begin{align*}
		M_{1}(\mathcal{NC}(G))&=(2mn-n)(2mn-n-1)^{2}-4(2mn-n-1)\dfrac{(mn-n)(mn-n-1)+mn(n-1)}{2} \\
		&~~~~~~~~~~~~~~~~~~~~~~~~~~~~~~~~~~~~~~~~~~~~~~~~~~+n(m-1)(mn-n-1)^{2}+mn(n-1)^{2} \\
		&= 5m^{3}n^{3}-9m^{2}n^{3}+4mn^{3} \\
		&= n^{3}(5m^{3}-9m^{2}+4m)
	\end{align*} 
	and 
	\begin{align*}
		M_{2}(\mathcal{NC}(G))&=\dfrac{(2mn-n)(2mn-n-1)^{3}}{2} + 2 \times \dfrac{((mn-n)(mn-n-1)+(mn^{2}-mn))^{2}}{4} \\
		&~~~~~~~~~~~~~~~~~~~~~~~~~~~~~~~~~~~~~-3 \times  \dfrac{(mn-n)(mn-n-1)+(mn^{2}-mn)}{2}(2mn-n-1)^{2} \\
		&~~~~~~~~~~~~~~~~~~~~~~~~~~~~~~~~~~~~~~~~~~+(2mn-n-\dfrac{3}{2})((mn-n)(mn-n-1)^{2}+mn(n-1)^{2}) \\
		&~~~~~~~~~~~~~~~~~~~~~~~~~~~~~~~~~~~~~~~~~~~~~~~-\dfrac{(mn-n)(mn-n-1)^{3}+mn(n-1)^{3}}{2} \\
		&= \dfrac{1}{2}(8m^{4}n^{4}-20m^{3}n^{4}+16m^{2}n^{4}-4mn^{4}) \\
		&= n^{4}(4m^{4}-10m^{3}+8m^{2}-2m).
	\end{align*}
	Also, $ |v(\mathcal{NC}(G) )|=2mn-n $ and $ |e( \mathcal{NC}(G))| = \binom{2mn-n}{2}-|e(\mathcal{C}(G))|=\frac{3m^{2}n^{2}-3mn^{2}}{2} $. We have 
	\[
	\dfrac{M_{1}(\mathcal{NC}(G))}{|v( \mathcal{NC}(G) )|} = \dfrac{mn^{2}(5m^{2}-9m+4)}{2m-1} \quad \text{ and } \quad \dfrac{M_{2}(\mathcal{NC}(G))}{|e(\mathcal{NC}(G))|} =  \dfrac{4n^{2}(2m^{3}-5m^{2}+4m-1)}{3(m -1)}.
	\]
	As such 
	\[
	\dfrac{M_{2}(\mathcal{NC}(G))}{|e( \mathcal{NC}(G) )|} - \dfrac{M_{1}(\mathcal{NC}(G))}{|v( \mathcal{NC}(G) )|} = \dfrac{n^{2}(m^{3}(m-6)+m(13m-12)+4))}{3(m-1)(2m-1)}:= \dfrac{n^2h(m)}{g(m)}.
	\]
	We have $g(m) > 0$ for all $m \geq 3$ and $h(m) > 0 $ for all $m \geq 6$. Also, $h(3) = 4 > 0$, $h(4)= 36$ and $h(5)= 144$.  Therefore, $\frac{n^2h(m)}{g(m)} > 0$ and so $\frac{M_{2}(\mathcal{NC}(G))}{|e(\mathcal{NC}(G) )|}  > \frac{M_{1}(\mathcal{NC}(G))}{|v( \mathcal{NC}(G) )|}$.
\end{proof}

\begin{cor}
	If $G=U_{6n} = \langle a, b: a^{2n}=b^{3}=1, a^{-1}ba=b^{-1}\rangle $, then 
	\[
	M_{1}(\mathcal{C}(G))= 2n(2n-1)^{2} + 3n(n-1)^{2}, M_{2}(\mathcal{C}(G))= \frac{2n(2n-1)^{3}+3n(n-1)^{3}}{2}, M_{1}(\mathcal{NC}(G))=66n^{3} 
	\]
	$ \text{ and } M_{2}(\mathcal{NC}(G))=120n^{4}. $ Further, $\frac{M_{2}(\Gamma(G))}{|e(\Gamma(G))|} > \frac{M_{1}(\Gamma(G))}{|v(\Gamma(G))|}$, where $\Gamma(G) = \mathcal{C}(G)$ or $\mathcal{NC}(G)$.
\end{cor}
\begin{proof}
	Since $\frac{U_{6n}}{Z(U_{6n})} \cong D_{6}$,  the result follows from Theorem \ref{thm-G/Z(G)=D_2m} considering $m=3.$
\end{proof}

\begin{cor}
	If $G=M_{2mn} = \langle a, b: a^{m}=b^{2n}=1, bab^{-1}=a^{-1} \rangle$ ($m \geq 3$ but not equal to $4$), then
	\[
	M_{1}(\mathcal{C}(G)) = \begin{cases}
		n(m-1)(mn-n-1)^{2} + mn(n-1)^{2}, & \text{when $m$ is odd} \\
		n(m-2)(mn-2n-1)^{2} + mn(2n-1)^{2}, & \text{when $m$ is even,}
	\end{cases} 
	\]
	\[
	M_{2}(\mathcal{C}(G))= \begin{cases}
		\dfrac{(mn-n)(mn-n-1)^{3}+mn(n-1)^{3}}{2}, & \text{when $m$ is odd} \\
		\dfrac{(mn-2n)(mn-2n-1)^{3}+mn(2n-1)^{3}}{2}, & \text{when $m$ is even,} 
	\end{cases}
	\]
	\[ 
	M_{1}(\mathcal{NC}(G)) = \begin{cases}
		n^{3}(5m^{3}-9m^{2}+4m), &\text{when m is odd} \\
		n^{3}(5m^{3}-18m^{2}+16m), &\text{when m is even}
	\end{cases}
	\]  
	and
	\[ 
	M_{2}(\mathcal{NC}(G))= \begin{cases}
		n^{4}(4m^{4}-10m^{3}+8m^{2}-2m), & \text{when m is odd} \\
		4n^{4}(m^{4}-5m^{3}+8m^{2}-4m), & \text{when m is even.} 
	\end{cases} 
	\]
Further, $\frac{M_{2}(\Gamma(G))}{|e(\Gamma(G))|} > \frac{M_{1}(\Gamma(G))}{|v(\Gamma(G))|}$, where $\Gamma(G) = \mathcal{C}(G)$ or $\mathcal{NC}(G)$.
\end{cor}
\begin{proof}
If $m$ is odd then $|Z(M_{2mn})| = n$ and $\frac{M_{2mn}}{Z(M_{2mn})} \cong D_{2m}$. Therefore, by  Theorem \ref{thm-G/Z(G)=D_2m}, we get 

$M_{1}(\mathcal{C}(G)) = 
	n(m-1)(mn-n-1)^{2} + mn(n-1)^{2}, M_{2}(\mathcal{C}(G))= \frac{(mn-n)(mn-n-1)^{3}+mn(n-1)^{3}}{2},$
	
$M_{1}(\mathcal{NC}(G)) = n^{3}(5m^{3}-9m^{2}+4m)$ and $M_{2}(\mathcal{NC}(G))= n^{4}(4m^{4}-10m^{3}+8m^{2}-2m)$.

\noindent Also, $\frac{M_{2}(\Gamma(G))}{|e(\Gamma(G))|} > \frac{M_{1}(\Gamma(G))}{|v(\Gamma(G))|}$, where $\Gamma(G) = \mathcal{C}(G)$ or $\mathcal{NC}(G)$.

If $m$ is even the $|Z(M_{2mn})| = 2n$ and $\frac{M_{2mn}}{Z(M_{2mn})} \cong D_{2\times \frac{m}{2}}$. Therefore, putting $n = 2n$ and $m = \frac{m}{2}$ in  Theorem \ref{thm-G/Z(G)=D_2m}, we get 

$M_{1}(\mathcal{C}(G)) = n(m-2)(mn-2n-1)^{2} + mn(2n-1)^{2}, M_{2}(\mathcal{C}(G))= \frac{(mn-2n)(mn-2n-1)^{3}+mn(2n-1)^{3}}{2}$,

$M_{1}(\mathcal{NC}(G)) = n^{3}(5m^{3}-18m^{2}+16m), M_{2}(\mathcal{NC}(G))= 4n^{4}(m^{4}-5m^{3}+8m^{2}-4m)$.

\noindent Also, $\frac{M_{2}(\Gamma(G))}{|e(\Gamma(G))|} > \frac{M_{1}(\Gamma(G))}{|v(\Gamma(G))|}$, where $\Gamma(G) = \mathcal{C}(G)$ or $\mathcal{NC}(G)$.
\end{proof}

\begin{theorem}\label{thm-ZpXZp}
	Let $G$ be a finite group such that $ \frac{G}{Z(G)} \cong \mathbb{Z}_{p} \times \mathbb{Z}_{p}$, where $p$ is a prime. Then $M_{1}(\mathcal{C}(G))=(pn-n)(p+1)(pn-n-1)^{2}$, $M_{2}(\mathcal{C}(G))=\frac{1}{2}(p+1)(pn-n)(pn-n-1)^{3}$, $M_{1}(\mathcal{NC}(G))$ $=(p+1)(pn-n)(p^{4}n^{2}-2p^{3}n^{2}+p^{2}n^{2}) $ and 
	$M_{2}(\mathcal{NC}(G))=\frac{1}{2}(p+1)^{2}(pn-n)^{2}(p^{4}n^{2}-2p^{3}n^{2}+p^{2}n^{2})$, where $n = |Z(G)|$. Further,
	$\frac{M_{2}(\Gamma(G))}{|e(\Gamma(G))|} = \frac{M_{1}(\Gamma(G))}{|v(\Gamma(G))|}$, where $\Gamma(G) = \mathcal{C}(G)$ or $\mathcal{NC}(G)$. 
\end{theorem} 
\begin{proof}
	It is well-known that $\mathcal{C}(G)=(p+1)K_{(p-1)n} $, where $ n=|Z(G)|$. As such, $ |v(\mathcal{C}(G))| = (p+1)(p-1)n=n(p^{2}-1) $ and $ |e(\mathcal{C}(G))|= (p+1) \cdot \binom{pn-n}{2}=\frac{(p+1)(pn-n)(pn-n-1)}{2}. $ Therefore, using \ref{thm1}, we get
	\[
	M_{1}(\mathcal{C}(G))=(p+1)(pn-n)(pn-n-1)^{2} \quad \text{ and } \quad M_{2}(\mathcal{C}(G)) = \dfrac{(p+1)(pn-n)(pn-n-1)^{3}}{2}.
	\]
	Also,
	\[
	\dfrac{M_{1}(\mathcal{C}(G))}{|v(\mathcal{C}(G))|}= (pn-n-1)^{2} = \dfrac{M_{2}(\mathcal{C}(G))}{|e(\mathcal{C}(G))|}.
	\]
	Using Theorem \ref{using commuting finding complement} we have 
	\begin{align*}
		M_{1}(\mathcal{NC}(G))&=(p^{2}n-n)(p^{2}n-n-1)^{2}-4(p^{2}n-n-1)\dfrac{(p^{2}n-n)(pn-n-1)}{2} \\
		&~~~~~~~~~~~~~~~~~~~~~~~~~~~~~~~~~~~~~~~~~~~~~~~~~~~~~~~~~~~~~~~+(p^{2}n-n)(pn-n-1)^{2} \\
		&= (p^{2}n-n)(p^{4}n^{2}-2p^{3}n^{2}+p^{2}n^{2}) \\
		&= (p+1)(pn-n)(p^{4}n^{2}-2p^{3}n^{2}+p^{2}n^{2})
	\end{align*} 
	and 
	\begin{align*}
		M_{2}(\mathcal{NC}(G))&=\dfrac{(p^{2}n-n)(p^{2}n-n-1)^{3}}{2}+2\dfrac{(p^{2}n-n)^{2}(pn-n-1)^{2}}{4} \\
		&~~~~~~~~~~~~~~~~~~~~~~~~~~~~~~~~~~~~~-3 \dfrac{(p^{2}n-n)(pn-n-1)}{2}(p^{2}n-n-1)^{2} \\
		&~~~~~~~~~~~~~~~~~~~~~~~~~~~~~~~~~~~~~~~~~~+(p^{2}n-n-\dfrac{3}{2})(p^{2}n-n)(pn-n-1)^{2} \\
		&~~~~~~~~~~~~~~~~~~~~~~~~~~~~~~~~~~~~~~~~~~~~~~~-\dfrac{(p^{2}n-n)(pn-n-1)^{3}}{2} \\
		&= \dfrac{p^{2}n-n}{2}(p^{6}n^{3}-3p^{5}n^{3}+3p^{4}n^{3}-p^{3}n^{3}) \\
		&=\dfrac{(p+1)^{2}(pn-n)^{2}(p^{4}n^{2}-2p^{3}n^{2}+p^{2}n^{2})}{2}.
	\end{align*}
	Also, $ |v(\mathcal{NC}(G) )|=p^{2}n-n $ and $ |e( \mathcal{NC}(G) )| = \binom{p^{2}n-n}{2}-|e(\mathcal{C}(G))|=\frac{(p^{2}n-n)(p^{2}n-pn)}{2} $. We have 
	\[
	\dfrac{M_{1}(\mathcal{NC}(G))}{|v( \mathcal{NC}(G) )|} = \dfrac{(p+1)(pn-n)(p^{4}n^{2}-2p^{3}n^{2}+p^{2}n^{2})}{p^{2}n-n}
	\]
	and
	\[
	\dfrac{M_{2}(\mathcal{NC}(G))}{|e(\mathcal{NC}(G))|} =  \dfrac{(p^{2}n-n)(p^{2}n-pn)(p^{4}n^{2}-2p^{3}n^{2}+p^{2}n^{2})}{(p^{2}n-n)(p^{2}n-pn)}.
	\]
	As such
	$$ 
	\dfrac{M_{1}(\mathcal{NC}(G))}{|v( \mathcal{NC}(G) )|} = p^{4}n^{2}-2p^{3}n^{2}+p^{2}n^{2} = \dfrac{M_{2}(\mathcal{NC}(G))}{|e( \mathcal{NC}(G))|}.
	$$ 
\end{proof}

\begin{theorem}
	Let $G$ be a finite group and $ \frac{G}{Z(G)} \cong Sz(2), $ where $Sz(2)$ is the Suzuki group presented by $\langle a, b : a^{5}=b^{4}=1, b^{-1}ab=a^{2}\rangle$. Then 
	$ M_{1}(\mathcal{C}(G))=4n(4n-1)^{2}+15n(3n-1)^{2}$, $M_{2}(\mathcal{C}(G))=\frac{1}{2}[4n(4n-1)^{3}+15n(3n-1)^{3}]$, $M_{1}(\mathcal{NC}(G))=4740n^{3}$ and $M_{2}(\mathcal{NC}(G))=37440n^{4}$, where $n=|Z(G)|$.
	Further, 	$\frac{M_{2}(\Gamma(G))}{|e(\Gamma(G))|} > \frac{M_{1}(\Gamma(G))}{|v(\Gamma(G))|}$, where $\Gamma(G) = \mathcal{C}(G)$ or $\mathcal{NC}(G)$.
\end{theorem}
\begin{proof}
	It is well-known that $\mathcal{C}(G) = K_{4n} \sqcup 5K_{3n}$, where $n=|Z(G)|$. As such,
	$ |v(\mathcal{C}(G))| =4n+5\cdot3n=19n $ and 
	$ |e(\mathcal{C}(G))| = \binom{4n}{2} + 5 \cdot\binom{3n}{2} = \frac{4n(4n-1)}{2}+5\cdot\frac{3n(3n-1)}{2}= \frac{4n(4n-1)+15n(3n-1)}{2}. $ 
	Therefore, using Theorem \ref{thm1}, we get
	\[
	M_{1}(\mathcal{C}(G))=4n(4n-1)^{2}+5\cdot3n(3n-1)^{2}=4n(4n-1)^{2}+15n(3n-1)^{2}
	\] 
	and
	\[
	M_{2}(\mathcal{C}(G))= \dfrac{4n(4n-1)^{3}}{2}+ 5 \cdot \dfrac{3n(3n-1)^{3}}{2}
	= \dfrac{4n(4n-1)^{3}+15n(3n-1)^{3}}{2}.
	\] 
	Also,
	\[ 
	\dfrac{M_{1}(\mathcal{C}(G))}{|v(\mathcal{C}(G))| }= \dfrac{4(4n-1)^{2}+15(3n-1)^{2}}{19} \quad \text{ and } \quad  \dfrac{M_{2}(\mathcal{C}(G))}{|e(\mathcal{C}(G))| }= \dfrac{4(4n-1)^{3}+15(3n-1)^{3}}{4(4n-1)+15(3n-1)}.
	\]
	We have $(7n-2)^{2}-4(4n-1)(3n-1) = n^2 > 0$. Therefore,
	\[
	(7n-2)^{2}-3(4n-1)(3n-1) > (4n-1)(3n-1).
	\]
	Multiplying both sides by $(7n - 2)$ we get
	\[
	(7n-2)^{3}-3(4n-1)(3n-1)(7n - 2) > (4n-1)(3n-1)(7n - 2).
	\]
	We have $(7n-2)^{3}-3(4n-1)(3n-1)(7n - 2) = (4n-1)^{3}+(3n-1)^{3}$ and so
	\[
	(4n-1)^{3}+(3n-1)^{3} > (4n-1)(3n-1)(7n - 2).
	\]
	Thus
	\[
	60(4n-1)^{3} + 60(3n-1)^{3} > 60(4n-1)(3n-1)(7n - 2).
	\]
	Again,
	\[
	60(4n-1)^{3} + 60(3n-1)^{3} = 76(4n-1)^{3} - 16(4n-1)^{3} + 285(3n-1)^{3} - 225(3n-1)^{3}
	\]
	and so
	\begin{align*}
		76(4n-1)^{3}   + 285(3n-1)^{3} & > 16(4n-1)^{3} + 225(3n-1)^{3} + 60(4n-1)(3n-1)(7n - 2)\\
		& = (4(4n-1)+15(3n-1))(4(4n-1)^{2}+15(3n-1)^{2}). 
	\end{align*}
	Therefore,
	\[
	\dfrac{4(4n-1)^{3}+15(3n-1)^{3}}{4(4n-1)+15(3n-1)} > \dfrac{4(4n-1)^{2}+15(3n-1)^{2}}{19}
	\]
and so  $\frac{M_{2}(\mathcal{C}(G))}{|e(\mathcal{C}(G))|} > \frac{M_{1}(\mathcal{C}(G))}{|v(\mathcal{C}(G))|}$.		
	
	Using Theorem \ref{using commuting finding complement} we have
	\begin{align*}
		M_{1}(\mathcal{NC}(G))&=19n(19n-1)^{2}-4(19n-1)\dfrac{(61n^{2}-19n)}{2}+4n(4n-1)^{2}+15n(3n-1)^{2} \\
		&= 6859n^{3}-2318n^{3}+64n^{3}+135n^{3} = 4740n^{3}
	\end{align*} 
	and 
	\begin{align*}
		M_{2}(\mathcal{NC}(G))&=\dfrac{19n(19n-1)^{3}}{2}+2\times \dfrac{(61n^{2}-19n)^{2}}{4}-3\times \dfrac{(61n^{2}-19n)}{2}(19n-1)^{2} \\ 
		&~~~~~~~~~~~~~~+(19n-\frac{3}{2})(4n(4n-1)^{2}+15n(3n-1)^{2})-\dfrac{4n(4n-1)^{3}+15n(3n-1)^{3}}{2}\\
		&= \frac{1}{2} \times 74880n^{4} = 37440n^{4}.
	\end{align*}
	Also, $ |v(\mathcal{NC}(G))|=19n$ and $ |e( \mathcal{NC}(G) )| = \binom{19n}{2}-|e(\mathcal{C}(G))|=150n^{2} $. We have $\frac{M_{1}(\mathcal{NC}(G))}{|v( \mathcal{NC}(G) )|} = \frac{4740n^{3}}{19n}$ and $ \frac{M_{2}(\mathcal{NC}(G))}{|e( \mathcal{NC}(G) )|} =  \frac{37440n^{4}}{150n^{2}}$. Therefore,  $
	\frac{M_{2}(\mathcal{NC}(G))}{|e( \mathcal{NC}(G) )|} >  \frac{M_{1}(\mathcal{NC}(G))}{|v( \mathcal{NC}(G) )|}
	$ since $ \frac{3744n^{2}}{15} > \frac{4740n^{2}}{19}$.
\end{proof}
Since $Sz(2)$ has trivial center we have the following corollary.
\begin{cor}\label{cor-Sz2}
	If $G \cong Sz(2)$ then 	$\frac{M_{2}(\Gamma(G))}{|e(\Gamma(G))|} > \frac{M_{1}(\Gamma(G))}{|v(\Gamma(G))|}$, where $\Gamma(G) = \mathcal{C}(G)$ or $\mathcal{NC}(G)$.
\end{cor} 

\subsection{Zagreb indices of $\mathcal{C}(G)$ and $\mathcal{NC}(G)$ for more groups}
In this subsection, we compute Zagreb indices of $\mathcal{C}(G)$ and $\mathcal{NC}(G)$ for Hanaki groups, certain general linear groups and projective special linear groups. However, we begin with the non-abelian group of order $pq$.

\begin{theorem}\label{thm-pq}
    Let $G$ be a finite non-abelian group of order $ pq $ where $p$ and $q$ are primes with $ p|(q-1) $. Then 
\[
M_{1}(\mathcal{C}(G))= (q-1)(q-2)^{2}+q(p-1)(p-2)^{2}, \quad M_{2}(\mathcal{C}(G))=\dfrac{(q-1)(q-2)^{3}+q(p-1)(p-2)^{3}}{2},
\]
\[
M_{1}(\mathcal{NC}(G))= p^{3}q^{3}-2p^{2}q^{2}-pq^{3}-p^{3}q^{2}+pq^{2}-3q^{2}-3qp^{2}+2q+p^{3}q+q^{3}-4
\]
and
\begin{align*}
    M_{2}(\mathcal{NC}(G))&=\dfrac{1}{2}\Big{(}p^{4}q^{4}-7p^{3}q^{3}+41p^{2}q^{2}-51pq+3p^{4}q^{2}+13q^{2}-16p^{3}q^{2}+14pq^{2}+2p^{2}q^{3} \\ &~~~~~~~~~~~~~~~~~~~~~~~~~~~~~~~~~~~~~~~~~~~~~~~-16pq^{3} +8p^{2}q-9q+2pq^{4}+2p^{3}q+p^{4}q+18 \Big{)}
\end{align*}
    Further, $\frac{M_{2}(\Gamma(G))}{|e(\Gamma(G))|} > \frac{M_{1}(\Gamma(G))}{|v(\Gamma(G))|}$, where $\Gamma(G) = \mathcal{C}(G)$ or $\mathcal{NC}(G)$.
\end{theorem} 
\begin{proof}
    It is well-known that $ \mathcal{C}(G)=K_{q-1} \sqcup qK_{p-1}. $ As such, $|v(\mathcal{C}(G))|= pq-1$ and 
\begin{align*}
    |e(\Gamma_{G}\mathcal{C}(G))|=\binom{q-1}{2}+ q \cdot \binom{p-1}{2}  
    = \dfrac{(q-1)(q-2)+q(p-1)(p-2)}{2}.
\end{align*}
    Therefore, using Theorem \ref{thm1}, we get
\begin{align*}
    M_{1}(\mathcal{C}(G))&= (q-1)(q-1-1)^{2}+q(p-1)(p-1-1)^{2}
    \\
    &=(q-1)(q-2)^{2}+q(p-1)(p-2)^{2}
\end{align*}
and
\begin{align*}
    M_{2}(\mathcal{C}(G))&= \dfrac{(q-1)(q-1-1)^{3}}{2}+q \cdot \dfrac{(p-1)(p-1-1)^{3}}{2} \\
    &= \dfrac{(q-1)(q-2)^{3}+q(p-1)(p-2)^{3}}{2}.    
\end{align*}
Also,
\[
\dfrac{M_{1}(\mathcal{C}(G))}{|v(\mathcal{C}(G))|}=\dfrac{(q-1)(q-2)^{2}+q(p-1)(p-2)^{2}}{pq-1} \quad \text{ and }
\] 
\[ 
\dfrac{M_{2}(\mathcal{C}(G))}{|e(\mathcal{C}(G))|}= \dfrac{(q-1)(q-2)^{3}+q(p-1)(p-2)^{3}}{(q-1)(q-2)+q(p-1)(p-2)}.
\]
    We have $(p+q-4)^{2}-4(p-2)(q-2) = (p - q)^2 > 0$ and so
\[
(p+q-4)^{2}-3(p-2)(q-2) > (p-2)(q-2).
\]
    Multiplying both sides by $(p + q - 4)$ we get
\[
(p+q-4)^{3}-3(p-2)(q-2)(p+q-4) > (p-2)(q-2)(p+q-4).
\]
    We have $(p+q-4)^{3}-3(p-2)(q-2)(p+q-4) = (q - 2)^3 + (p - 2)^3$ and so
\[
(q - 2)^3 + (p - 2)^3 > (p-2)(q-2)(p+q-4).
\]
    Multiplying both sides by $q(q - 1)(p - 1)$ we get
\[
f(p, q) := q(q - 1)(p - 1)(q - 2)^3 + q(q - 1)(p - 1)(p - 2)^3 > q(q - 1)(p - 1)(p-2)(q-2)(p+q-4).
\]
    Again,
\[
f(p, q) = (q - 1)(q - 2)^3(pq - 1) - (q - 1)^2(q - 2)^3 + (p - 1)(p-2)^3q(pq - 1) - (p - 1)^2(p - 2)^3q^2.
\]
Therefore,
\begin{align*}
    (q - 1)(q - 2)^3&(pq - 1)  + (p - 1)(p-2)^3q(pq - 1)\\
    > & (q - 1)^2(q - 2)^3 + (p - 1)^2(p - 2)^3q^2 + q(q - 1)(p - 1)(p-2)(q-2)(p+q-4)\\
    = & \Big{(} (q-1)(q-2)+q(p-1)(p-2)\Big{)}\Big{(}(q-1)(q-2)^{2}+q(p-1)(p-2)^{2}\Big{)}
\end{align*}
    and so
\[
\dfrac{(q-1)(q-2)^{3}+q(p-1)(p-2)^{3}}{(q-1)(q-2)+q(p-1)(p-2)} > \dfrac{(q-1)(q-2)^{2}+q(p-1)(p-2)^{2}}{pq-1}.
\]
Thus $\frac{M_{2}(\mathcal{C}(G))}{|e(\mathcal{C}(G))|} \geq \frac{M_{1}(\mathcal{C}(G))}{|v(\mathcal{C}(G))|}$.

    Using Theorem \ref{using commuting finding complement} we have
\begin{align*}
    M_{1}(\mathcal{NC}(G))&=(pq-1)(pq-2)^{2}-4(pq-2)\dfrac{(q-1)(q-2)+q(p-1)(p-2)}{2} \\
    &~~~~~~~~~~~~~~~~~~~~~~~~~~~~~~~~~~~~~~~~~~~~~~~~~~~~~~+(q-1)(q-2)^{2}+q(p-1)(p-2)^{2} \\
    &= p^{3}q^{3}-2p^{2}q^{2}-pq^{3}-p^{3}q^{2}+pq^{2}-3q^{2}-3qp^{2}+2q+p^{3}q+q^{3}-4
\end{align*} 
and 
\begin{align*}
    M_{2}(\mathcal{NC}(G))&=\dfrac{(pq-1)(pq-2)^{3}}{2}+2\dfrac{[(q-1)(q-2)+q(p-1)(p-2)]^{2}}{4} \\
    &~~~~~~~~~~~~~~~~~~~~~~~~~~~~~~~~~~~~~~~~-3\dfrac{(q-1)(q-2)+q(p-1)(p-2)}{2}(pq-2)^{2} \\
    &~~~~~~~~~~~~~~~~~~~~~~~~~~~~~~~~~~~~~~~~+(pq-1-\dfrac{3}{2})[(q-1)(q-2)^{2}+q(p-1)(p-2)^{2}]\\ 
    &~~~~~~~~~~~~~~~~~~~~~~~~~~~~~~~~~~~~~~~~-\dfrac{(q-1)(q-2)^{3}+q(p-1)(p-2)^{3}}{2} \\
    &= \dfrac{1}{2}\Big{(}p^{4}q^{4}-7p^{3}q^{3}+41p^{2}q^{2}-51pq+3p^{4}q^{2}+13q^{2}-16p^{3}q^{2}+14pq^{2}+2p^{2}q^{3} \\ &~~~~~~~~~~~~~~~~~~~~~~~~~~~~~~~~~~~~~~~~~~~~~~~~-16pq^{3}+8p^{2}q-9q+2pq^{4}+2p^{3}q+p^{4}q+18 \Big{)}
\end{align*}
    Also, $ |v(\mathcal{NC}(G) )|=pq-1 $ and $ |e( \mathcal{NC}(G) )| = \binom{pq-1}{2}-|e(\mathcal{C}(G))|=\frac{p^{2}q^{2}-p^{2}q-q^{2}+q}{2} $. 
As such,
\begin{align*}
    &\dfrac{M_{2}(\mathcal{NC}(G))}{|e( \mathcal{NC}(G) )|} - \dfrac{M_{1}(\mathcal{NC}(G))}{|v( \mathcal{NC}(G) )|} \\
    &~~~=\dfrac{2p^{4}q^{4}(p-3)+p^{2}q^{4}(pq-2p-14)+p^{3}q^{3}(q^{2}-15p)+p^{2}q(q^{4}+p^{3}q^{2}-p^{2})+69pq}{pq^{2}(p^{2}q-q-p-p^{2})+q(q-1)+pq(p+q)} \\
    &~~~+ \dfrac{pq^{2}(6pq-23)+p^{4}q^{2}(2p-4)+p^{2}q^{2}(51pq-83)+p^{3}q(23q-2)+pq(28q^{2}-12p)}{pq^{2}(p^{2}q-q-p-p^{2})+q(q-1)+pq(p+q)} \\
    &~~~=\dfrac{A(p, q)}{B(p, q)},
\end{align*}
where $A(p, q) := 2p^{4}q^{4}(p-3)+p^{2}q^{4}(pq-2p-14)+p^{3}q^{3}(q^{2}-15p)+p^{2}q(q^{4}+p^{3}q^{2}-p^{2})+69pq$ $+ pq^{2}(6pq-23)+p^{4}q^{2}(2p-4)+p^{2}q^{2}(51pq-83)+p^{3}q(23q-2)+pq(28q^{2}-12p)$ and $B(p, q) := pq^{2}(p^{2}q-q-p-p^{2})+q(q-1)+pq(p+q) = pq^{2}(q(p^2 - 1) - p(p + 1)) + q(q-1)+pq(p+q)$. Since $q(p^{2}-1)>p(p+1)$ and $q > 1$ we have $B(p, q)>0$. In order to determine whether $A(p, q) > 0$ or not we consider the following cases.

\noindent \textbf{Case 1.} $p=2$
    
    We have $A(2, q)= q^{4}(20q-104)+q^{2}(280q-194)+58q$ and so  $A(2, q)>0$ for $q \geq 7$. Also  $A(2, 3)=2424$ and $A(2, 5)=27940$. 

    \noindent \textbf{Case 2.} $p \geq 3$.

    We have $p - 3 \geq 0$, $p(q-2)>14$, $q^{2}>15p$, $q^{2}(q^{2}+p^{3})>p^{2}$, $6pq > 23$, $2p > 4$, $51pq - 83$, $23q > 2$ and $28q^2 > 12p$ and so  $A(p, q)>0$.

     Therefore, in all the case, $A(p, q) > 0$ and hence $\frac{A(p, q)}{B(p, q)} > 0$. That is, $\frac{M_{2}(\mathcal{NC}(G))}{|e( \mathcal{NC}(G) )|}  \geq \frac{M_{1}(\mathcal{NC}(G))}{|v( \mathcal{NC}(G))|}$. 
\end{proof}
\begin{theorem}
	Let $F=GF(2^{n}), n \geq 2$ and $\nu$ be the Frobenius automorphism of $F$, i.e. $\nu(x)= x^{2}~~ \forall x \in F.$ Then the first and second Zagreb indices of the commuting and non-commuting graph of the group \[A(n,\nu) = \Bigg\lbrace U(a, b)=
	\begin{bmatrix}
		1 & 0 & 0 \\
		a & 1 & 0 \\
		b & \nu(a) & 1 
	\end{bmatrix} : a, b \in F \Bigg\rbrace \] are given by $M_{1}(\mathcal{C}(A(n,\nu)))=2^{n}(2^{n}-1)^{3}$, $M_{2}(\mathcal{C}(A(n,\nu)))=2^{n-1}(2^{n}-1)^{4}$, $M_{1}(\mathcal{NC}(A(n,\nu)))=2^{5n}(2^{n}-5)+2^{3n+2}(2^{n+1}-1)$ and $M_{2}(\mathcal{NC}(A(n,\nu)))=2^{7n}(2^{n-1}-3)-2^{6n}(2^{n-1}-9)-2^{4n+1}(5\cdot 2^{n}-2) $.
	Further, $\frac{M_{2}(\Gamma(A(n,\nu)))}{|e(\Gamma(A(n,\nu)))|} = \frac{M_{1}(\Gamma(A(n,\nu)))}{|v(\Gamma(A(n,\nu)))|}$, where $\Gamma(A(n,\nu)) = \mathcal{C}(A(n,\nu))$ or $\mathcal{NC}(A(n,\nu))$.
\end{theorem}
\begin{proof} 
    It is well-known that $\mathcal{C}(A(n,\nu))= (2^{n}-1)K_{2^{n}} $. As such, $ |v(\mathcal{C}(A(n,\nu)))|=(2^{n}-1)2^{n}=2^{2n}-2^{n} $ and $ |e(\mathcal{C}(A(n,\nu)))|=(2^{n}-1)\binom{2^{n}}{2}=2^{n-1}(2^{n}-1)^{2}$. Therefore, using Theorem \ref{thm1}, we get
\[
M_{1}(\mathcal{C}(A(n,\nu)))=(2^{n}-1)2^{n}(2^{n}-1)^{2}=2^{n}(2^{n}-1)^{3} \quad \text{ and}
\] 
\[
M_{2}(\mathcal{C}(A(n,\nu)))=(2^{n}-1) \times \dfrac{2^{n}(2^{n}-1)^{3}}{2}=2^{n-1}(2^{n}-1)^{4}.
\]
Therefore,
\[
\dfrac{M_{1}(\mathcal{C}(A(n,\nu)))}{|v(\mathcal{C}(A(n,\nu)))|}=	(2^{n}-1)^{2} = \dfrac{M_{2}(\mathcal{C}(A(n,\nu)))}{|e(\mathcal{C}(A(n,\nu)))|}.
\]
    Using Theorem \ref{using commuting finding complement} we have
\begin{align*}
    M_{1}(\mathcal{NC}(A(n,\nu)))&=(2^{2n}-2^{n})(2^{2n}-2^{n}-1)^{2}-4(2^{2n}-2^{n}-1)(2^{n-1}(2^{n}-1)^{2})+2^{n}(2^{n}-1)^{3}\\
    &= 2^{6n}-5\cdot2^{5n}+8\cdot2^{4n}-4\cdot2^{3n} \\
    &= 2^{5n}(2^{n}-5)+2^{3n+2}(2^{n+1}-1)
 \end{align*} 
    and 
 \begin{align*}
    M_{2}(\mathcal{NC}(A(n,\nu)))&=\dfrac{(2^{2n}-2^{n})(2^{2n}-2^{n}-1)^{3}}{2}+2\cdot2^{2n-2}(2^{n}-1)^{4} \\
    &~~~-3\cdot2^{n-1}(2^{n}-1)^{2}(2^{2n}-2^{n}-1)^{2} +\left(2^{2n}-2^{n}-\dfrac{3}{2}\right)2^{n}(2^{n}-1)^{3}-\dfrac{2^{n}(2^{n}-1)^{4}}{2} \\
    &= 2^{8n-1}-3\cdot2^{7n}-2^{7n-1}+9\cdot2^{6n}-10\cdot2^{5n}+2^{4n+2} \\
    &=2^{7n}(2^{n-1}-3)-2^{6n}(2^{n-1}-9)-2^{4n+1}(5\cdot 2^{n}-2).
\end{align*}
    Also, $ |v(\mathcal{NC}(A(n,\nu)))|=2^{2n}-2^{n} $ and $ |e(\mathcal{NC}(A(n,\nu)) )| = \binom{2^{2n}-2^{n}}{2}-|e(\mathcal{C}(A(n,\nu)))|=2^{n}(2^{n}-2)(2^{2n}-2^{n})$. Therefore 
 \begin{align*}
\dfrac{M_{2}(\mathcal{NC}(A(n,\nu)))}{|e( \mathcal{NC}(A(n,\nu)) )|} & = \dfrac{2^{7n}(2^{n-1}-3)-2^{6n}(2^{n-1}-9)-2^{4n+1}(5\cdot 2^{n}-2)}{2^{n}(2^{n}-2)(2^{2n}-2^{n})}\\
& = \dfrac{2^{5n}(2^{n}-5)+2^{3n+2}(2^{n+1}-1)}{2^{2n}-2^{n}}\\
& = \dfrac{M_{1}(\mathcal{NC}(A(n,\nu)))}{|v(\mathcal{NC}(A(n,\nu)))|}.
\end{align*}  
\end{proof}
\begin{theorem}
    Let $F=GF(p^{n}), p$ be a prime. Then the first and second Zagreb indices of the commuting and non-commuting graph of the group 
\[
A(n, p) = \Bigg\lbrace v(a, b, c)=
	\begin{bmatrix}
		1 & 0 & 0 \\
		a & 1 & 0 \\
		b & c & 1 
	\end{bmatrix} : a, b, c \in F \Bigg\rbrace 
\] 
    are given by 
\[
M_{1}(\mathcal{C}(A(n, p)))=p^{n}(p^{2n}-1)(p^{2n}-p^{n}-1)^{2}, 
M_{2}(\mathcal{C}(A(n, p)))=\dfrac{p^{n}(p^{2n}-1)(p^{2n}-p^{n}-1)^{3}}{2},
\]
\[
M_{1}(\mathcal{NC}(A(n, p)))=p^{8n}(p^{n}-2)+p^{5n}(2p^{n}-1) 
\]
and
\[
M_{2}(\mathcal{NC}(A(n, p)))=\dfrac{(p^{3n}-p^{n})[p^{8n}(p^{n}-3)+p^{6n}(3p^{n}-1)]}{2}.
\]
    Further, $\frac{M_{2}(\Gamma(A(n, p)))}{|e(\Gamma(A(n, p)))|} = \frac{M_{1}(\Gamma(A(n, p)))}{|v(\Gamma(A(n, p)))|}$, where $\Gamma(A(n, p)) = \mathcal{C}(A(n, p))$ or $\mathcal{NC}(A(n, p))$.
\end{theorem}
\begin{proof}
	It is well-known that $\mathcal{C}(A(n, p))=(p^{n}+1)K_{p^{2n}-p^{n}}.$
	As such, $|v(\mathcal{C}(A(n, p)))|=(p^{n}+1)(p^{2n}-p^{n})=
p^{3n}-p^{n}$ and 
$|e(\mathcal{C}(A(n, p)))|=(p^{n}+1)\binom{p^{2n}-p^{n}}{2}= 
\frac{p^{n}(p^{2n}-1)(p^{2n}-p^{n}-1)}{2}.$
    Therefore, using Theorem \ref{thm1}, we get
\[
M_{1}(\mathcal{C}(A(n, p)))=(p^{n}+1)(p^{2n}-p^{n})(p^{2n}-p^{n}-1)^{2}=p^{n}(p^{2n}-1)(p^{2n}-p^{n}-1)^{2}
\]  
and
\[
M_{2}(\mathcal{C}(A(n, p)))=(p^{n}+1) \dfrac{(p^{2n}-p^{n})(p^{2n}-p^{n}-1)^{3}}{2}=\dfrac{p^{n}(p^{2n}-1)(p^{2n}-p^{n}-1)^{3}}{2}.
\]
    Also,
\[
\dfrac{M_{1}(\mathcal{C}(A(n, p)))}{|v(\mathcal{C}(A(n, p)))|}=
(p^{2n}-p^{n}-1)^{2} = \dfrac{M_{2}(\mathcal{C}(A(n, p)))}{|e(\mathcal{C}(A(n, p)))|}.
\]
    Using Theorem \ref{using commuting finding complement} we have
\begin{align*}
    M_{1}(\mathcal{NC}(A(n, p)))&=(p^{3n}-p^{n})(p^{3n}-p^{n}-1)^{2}-4(p^{3n}-p^{n}-1)\dfrac{p^{n}(p^{2n}-1)(p^{2n}-p^{n}-1)}{2}\\
    &~~~~~~~~~~~~~~~~~~~~~~~~~~~~~~~~~~~~~~~~~~~~~~~~~~~~~~~~~~~~~~~~~~~~+p^{n}(p^{2n}-1)(p^{2n}-p^{n}-1)^{2}\\
    &= p^{9n}-2p^{8n}-p^{5n}+2p^{6n} \\
    &= p^{8n}(p^{n}-2)+p^{5n}(2p^{n}-1)
\end{align*} 
    and 
\begin{align*}
    M_{2}(\mathcal{NC}(A(n, p)))&=\dfrac{(p^{3n}-p^{n})(p^{3n}-p^{n}-1)^{3}}{2}+2\dfrac{(p^{3n}-p^{n})^{2}(p^{2n}-p^{n}-1)^{2}}{4}\\
    &~~~~~~~~~~~~~~~~~~~~~~~~~~~~~~~~-3\dfrac{(p^{3n}-p^{n})(p^{2n}-p^{n}-1)}{2}(p^{3n}-p^{n}-1)^{2}\\
    &~~~~~~~~~~~~~~~~~~~~~~~~~~~~~~~~+(p^{3n}-p^{n}-\dfrac{3}{2})(p^{3n}-p^{n})(p^{2n}-p^{n}-1)^{2}\\
    &~~~~~~~~~~~~~~~~~~~~~~~~~~~~~~~~-\dfrac{(p^{3n}-p^{n})(p^{2n}-p^{n}-1)^{3}}{2} \\
    &= \dfrac{(p^{3n}-p^{n})}{2}(p^{9n}+3p^{7n}-3p^{8n}-p^{6n}) \\
    &=\dfrac{(p^{3n}-p^{n})(p^{8n}(p^{n}-3)+p^{6n}(3p^{n}-1))}{2}.
 \end{align*}
    Also, $ |v(\mathcal{NC}(A(n, p)))|=p^{3n}-p^{n} $ and $ |e( \mathcal{NC}(A(n, p)))| = \binom{p^{3n}-p^{n}}{2}-|e(\mathcal{C}(A(n, p)))|=\frac{p^{2n}}{2}(p^{n}-1)(p^{3n}-p^{n})$. Therefore
\begin{align*}
\dfrac{M_{2}(\mathcal{NC}(A(n, p)))}{|e(\mathcal{NC}(A(n, p)) )|} &=  \dfrac{(p^{3n}-p^{n})(p^{8n}(p^{n}-3)+p^{6n}(3p^{n}-1))}{p^{2n}(p^{n}-1)(p^{3n}-p^{n})}\\
& = \dfrac{p^{8n}(p^{n}-2)+p^{5n}(2p^{n}-1)}{p^{3n}-p^{n}}\\
& =  \dfrac{M_{1}(\mathcal{NC}(A(n, p)))}{|v( \mathcal{NC}(A(n, p)) )|}. 
\end{align*}    
\end{proof}
\begin{theorem}
    Let $G=GL(2, q)$ (the general linear group), where $q=p^{n}>2$ and $p$ is a prime integer. Then 
$$
M_{1}(\mathcal{C}(G))= q(q-1)(q^{6}-4q^{5}+4q^{4}+2q^{3}-4q^{2}+q-1),
$$
$$
M_{2}(\mathcal{C}(G))= \dfrac{q(q-1)}{2} (q^{8}-6q^{7}+14q^{6}-15q^{5}+3q^{4}+12q^{3}-16q^{2}+9q-1),
$$ 
$$ 
M_{1}(\mathcal{NC}(G))=(q-1)(q^{11}-2q^{10}-4q^{9}+9q^{8}+5q^{7}-15q^{6}+q^{5}+7q^{4}-2q^{3}+q^{2}-q) \quad \text{ and } \quad
$$ 
\begin{align*}
    M_{2}(\mathcal{NC}(G))&= \dfrac{q(q-1)}{2}(q^{14}-3q^{13}-4q^{12}+19q^{11}-47q^{9}+28q^{8}+43q^{7}-50q^{6}+11q^{5}+4q^{4} \\ &~~~~~~~~~~~~~~~~~~~~~~~~~~~~~~~~~~~~~~~~~~~~~~~~~~~~~~~~~~~~~~~~~~~~~~~~~~~~~~~~~~~~~~~-12q^{3}+19q^{2}-11q+2).
\end{align*}  
    Further, $\frac{M_{2}(\Gamma(G))}{|e(\Gamma(G))|} > \frac{M_{1}(\Gamma(G))}{|v(\Gamma(G))|}$, where $\Gamma(G) = \mathcal{C}(G)$ or $\mathcal{NC}(G)$.
\end{theorem}
\begin{proof}
    It is well-known that $ |G|=(q^{2}-1)(q^{2}-q), ~ |Z(G)|=q-1 $ and $ \mathcal{C}(G)= \frac{q(q+1)}{2}K_{q^{2}-3q+2} \sqcup \frac{q(q-1)}{2} K_{q^{2}-q} \sqcup (q+1) K_{q^{2}-2q+1} $. As such, $|v(\mathcal{C}(G))|=(q-1)(q^{3}-q-1)$ and $ |e(\mathcal{C}(G))|= \frac{q(q-1)}{2} ( q^{4} - 2q^{3} - q^{2} + 2q + 1 ).$ Therefore, using Theorem 2.1, we get 
\begin{align*}
    M_{1}(\mathcal{C}(G))&=\dfrac{q(q+1)}{2}(q^{2}-3q+2)(q^{2}-3q+1)^{2} + \dfrac{q(q-1)}{2}(q^{2}-q)(q^{2}-q-1)^{2} \\
    &~~~~~~~~~~~~~~~~~~~~~~~~~~~~~~~~~~~~~~~~~~~~~~~~~~~~~~~~~~~~~~+ (q+1)(q^{2}-2q+1)(q^{2}-2q)^{2}\\
    &= q(q-1)(q^{6}-4q^{5}+4q^{4}+2q^{3}-4q^{2}+q-1)
\end{align*}
and 
\begin{align*}
    M_{2}(\mathcal{C}(G))&= \dfrac{q(q+1)}{2}(q^{2}-3q+2)\dfrac{(q^{2}-3q+1)^{3}}{2} + \dfrac{q(q-1)}{2}(q^{2}-q)\dfrac{(q^{2}-q-1)^{3}}{2} \\ 
    &~~~~~~~~~~~~~~~~~~~~~~~~~~~~~~~~~~~~~~~~~~~~~~~~~~~~~~~~~~~~~~~~~~~~~~~~~~ + (q+1)(q^{2}-2q+1)\dfrac{(q^{2}-2q)^{3}}{2}\\
    &= \dfrac{q(q-1)}{2} (q^{8}-6q^{7}+14q^{6}-15q^{5}+3q^{4}+12q^{3}-16q^{2}+9q-1).
\end{align*}
We have 
$$ 
\dfrac{ M_{1}(\mathcal{C}(G))}{|v(\mathcal{C}(G))|}= \dfrac{q(q-1)(q^{6}-4q^{5}+4q^{4}+2q^{3}-4q^{2}+q-1)}{(q-1)(q^{3}-q-1)}
$$ 
and 
$$
\dfrac{M_{2}(\mathcal{C}(G))}{|e(\mathcal{C}(G))|}= \dfrac{q^{8}-6q^{7}+14q^{6}-15q^{5}+3q^{4}+12q^{3}-16q^{2}+9q-1}{q^{4} - 2q^{3} - q^{2} + 2q + 1} .
$$
Therefore, 
\begin{align*}
    \dfrac{M_{2}(\mathcal{C}(G))}{|e(\mathcal{C}(G))|} & - \dfrac{ M_{1}(\mathcal{C}(G))}{|v(\mathcal{C}(G))|} \\
    &=\dfrac{2q^{8}(q-5)+q^{5}(14q^{2}-13)+q^{3}(24q^{3}-q+4)+q(8q-7)+1}{q^{5}(q^{2}-2q-2)+q(3q^{3}-q-3)+(4q^{3}-1)} 
    := \dfrac{f(q)}{g(q)}.   
\end{align*}
Since $q > 2 $ we have $q - 2 \geq 1, q^{3}-q>1 $ and $ 4q^{3}-1>0. $ As such, $q(q-2)=q^{2}-2q>2$ and $ 3q^{3}-3=3(q^{3}-1)>q$ and so $ g(q)>0.$ For $ q>3 $ we have $ q-5>0, 14q^{2}-13>0, 24q^{3}-q+4>0 $ $ 8q-7>0  $ and so $ f(q)>0$. Also  $ f(3)=18,7880$. Therefore, $ \frac{f(q)}{g(q)}>0$. Thus, $\frac{M_{2}(\mathcal{C}(G))}{|e(\mathcal{C}(G))|} > \frac{M_{1}(\mathcal{C}(G))}{|v(\mathcal{C}(G))|}$.

    Using Theorem \ref{using commuting finding complement} we have
\begin{align*}
    M_{1}(\mathcal{NC}(G))&=(q-1)(q^{3}-q-1)((q-1)(q^{3}-q-1)-1)^{2} \\ 
    &~~~~-4((q-1)(q^{3}-q-1)-1)\dfrac{q(q-1)}{2} ( q^{4} - 2q^{3} - q^{2} + 2q + 1 )\\
    &~~~~~+ q(q-1)(q^{6}-4q^{5}+4q^{4}+2q^{3}-4q^{2}+q-1) \\
    &= (q-1)(q^{11}-2q^{10}-4q^{9}+9q^{8}+5q^{7}-15q^{6}+q^{5}+7q^{4}-2q^{3}+q^{2}-q)
\end{align*} 
and  
\begin{align*}
    M_{2}(\mathcal{NC}(G))&=\dfrac{(q-1)(q^{3}-q-1)[(q-1)(q^{3}-q-1)-1]^{3}}{2} \\
    &~~~~~~~~~~~~~~+2\times \dfrac{q^{2}(q-1)^{2}( q^{4}-2q^{3}-q^{2}+2q+1)}{4} \\ 
    &~~~~~~~~~~~~~~-3\times \dfrac{q(q-1)}{2}( q^{4}-2q^{3}-q^{2}+2q+1)((q-1)(q^{3}-q-1)-1)^{2} \\  
    &~~~~~~~~~~~~~+((q-1)(q^{3}-q-1)-\dfrac{3}{2})(q(q-1)(q^{6}-4q^{5}+4q^{4}+2q^{3}-4q^{2}+q-1)) \\
    &~~~~~~~~~~~~~~-\dfrac{q(q-1)}{2} (q^{8}-6q^{7}+14q^{6}-15q^{5}+3q^{4}+12q^{3}-16q^{2}+9q-1)\\
    &= \dfrac{q(q-1)}{2}(q^{14}-3q^{13}-4q^{12}+19q^{11}-47q^{9}+28q^{8}+43q^{7}-50q^{6}+11q^{5}+4q^{4} \\ &~~~~~~~~~~~~~~~~~~~~~~~~~~~~~~~~~~~~~~~~~~~~~~~~~~~~~~~~~~~~~~~~~~~~~~~~~~~~~~~~~~~~~~~~~-12q^{3}+19q^{2}-11q+2) \\
    &:=\frac{q}{2}A(q),
\end{align*}
where $ A(q)= (q-1)(q^{14}-3q^{13}-4q^{12}+19q^{11}-47q^{9}+28q^{8}+43q^{7}-50q^{6}+11q^{5}+4q^{4}-12q^{3}+19q^{2}-11q+2)$. Also, $ |v(\mathcal{NC}(G))|=(q-1)(q^{3}-q-1)$ and $ |e( \mathcal{NC}(G) )| = \binom{(q-1)(q^{3}-q-1)}{2}-|e(\mathcal{C}(G))|=\frac{q}{2}(q^{7}-2q^{6}-2q^{5}+5q^{4}-4q^{2}+q^{3}+1) $. We have 
 $$
 \dfrac{M_{1}(\mathcal{NC}(G))}{|v( \mathcal{NC}(G) )|} =\dfrac{(q^{11}-2q^{10}-4q^{9}+9q^{8}+5q^{7}-15q^{6}+q^{5}+7q^{4}-2q^{3}+q^{2}-q)}{(q^{3}-q-1)}
 $$ 
 and 
\begin{align*}
    \dfrac{M_{2}(\mathcal{NC}(G))}{|e( \mathcal{NC}(G) )|} &=  \dfrac{A(q)}{(q^{7}-2q^{6}-2q^{5}+5q^{4}-4q^{2}+q^{3}+1)}.
\end{align*}
As such 
 \begin{align*}
    &\dfrac{M_{2}(\mathcal{NC}(G))}{|e(\mathcal{NC}(G))|} - \dfrac{M_{1}(\mathcal{NC}(G))}{|v( \mathcal{NC}(G) )|} \\ &~~~~~=\dfrac{q^{11}(q-5)+q^{8}(14q-35)+q^{5}(15q-12)+q^{4}(5q^{6}-4)+q^{3}(18q^{4}-5)+16q^{2}-10q+2}{(q^{5}(q^{2}-2q-2)+q^{2}(5q^{2}-4)+q^{3}+1)(q^{3}-q-1)} \\ 
    &~~~~~:=\dfrac{f(q)}{g(q)}.
\end{align*} 
    We have $g(q) > 0$, $f(3)=33920$ and $f(4)=2767770$. For $q \geq 5$ we have $f(q) > 0 $. Therefore, $\frac{f(q)}{g(q)}>0.$ Thus, $\frac{M_{2}(\mathcal{NC}(G))}{|e(\mathcal{NC}(G))|} > \frac{M_{1}(\mathcal{NC}(G))}{|v( \mathcal{NC}(G) )|}$.  
\end{proof}
\begin{theorem}\label{thm-PSL}
    If $ G=PSL(2,2^{k})$ (the projective special linear group), where $k\geq 2$, then 
\[
M_{1}(\mathcal{C}(G))=2^{5k}-4 \cdot 2^{4k} + 4 \cdot 2^{3k} +4 \cdot 2^{2k} -5 \cdot 2^{k} -4,
\]  
\[
M_{2}(\mathcal{C}(G))=\dfrac{2^{6k}-6\cdot2^{5k}+14\cdot2^{4k}-9\cdot2^{3k}-15\cdot2^{2k}+15\cdot2^{k}+8}{2},
\]
\[
M_{1}(\mathcal{NC}(G))=2^{9k}-5\cdot2^{7k}-2^{6k}+9\cdot2^{5k}-5\cdot2^{3k}-3\cdot2^{2k}+3\cdot2^{k}
\]
and
\[
M_{2}(\mathcal{NC}(G))=\frac{1}{2}(2^{12k}-7\cdot2^{10k}-2^{9k}+21\cdot2^{8k}-26\cdot2^{6k}-2\cdot2^{5k}+15\cdot2^{4k}+3\cdot2^{3k}+6\cdot2^{2k}-8\cdot2^{k}).
\]
    Further, $\frac{M_{2}(\Gamma(G))}{|e(\Gamma(G))|} > \frac{M_{1}(\Gamma(G))}{|v(\Gamma(G))|}$, where $\Gamma(G) = \mathcal{C}(G)$ or $\mathcal{NC}(G)$.
\end{theorem}
\begin{proof}
    It is well-known that $ \mathcal{C}(G)=(2^{k}+1)K_{2^{k}-1} \sqcup 2^{k-1}(2^{k}+1)K_{2^{k}-2} \sqcup 2^{k-1}(2^{k}-1)K_{2^{k}}. $ As such, $ |v(\mathcal{C}(G))|=(2^{k}+1)(2^{k}-1)+2^{k-1}(2^{k}+1)(2^{k}-2)+2^{k-1}(2^{k}-1)2^{k}=2^{3k}-2^{k}-1 $ and
\begin{align*}
    |e(\mathcal{C}(G))|&= \dfrac{(2^{k}+1)(2^{k}-1)(2^{k}-2)}{2}+\dfrac{2^{k-1}(2^{k}+1)(2^{k}-2)(2^{k}-3)}{2}+\dfrac{2^{k-1}(2^{k}-1)2^{k}(2^{k}-1)}{2} \\
    &= \dfrac{2^{4k}-2\cdot2^{3k}-2\cdot2^{2k}+3\cdot2^{k}+2}{2}.
\end{align*}
    Therefore, using Theorem \ref{thm1}, we have
\begin{align*}
    M_{1}(\mathcal{C}(G))&=(2^{k}+1)(2^{k}-1)(2^{k}-1-1)^{2}+2^{k-1}(2^{k}+1)(2^{k}-2)(2^{k}-2-1)^{2} \\
    &~~~~~~~~~~~~~~~~~~~~~~~~~~~~~~~~~~~~~~~~~~~~~~~~~~~~~~~~~~~~~~~~~~~~~~~~~+2^{k-1}(2^{k}-1)2^{k}(2^{k}-1)^{2} \\
    &=(2^{k}+1)(2^{k}-1)(2^{k}-2)^{2}+2^{k-1}(2^{k}+1)(2^{k}-2)(2^{k}-3)^{2}+2^{k-1}(2^{k}-1)2^{k}(2^{k}-1)^{2} \\
    &= 2^{5k}-4 \cdot 2^{4k} + 4 \cdot 2^{3k} +4 \cdot 2^{2k} -5 \cdot 2^{k} -4
\end{align*}
and
\begin{align*}
    M_{2}(\mathcal{C}(G))&= (2^{k}+1) \dfrac{(2^{k}-1)(2^{k}-1-1)^{3}}{2}+2^{k-1}(2^{k}+1)\dfrac{(2^{k}-2)(2^{k}-2-1)^{3}}{2}\\
    &~~~~~~~~~~~~~~~~~~~~~~~~~~~~~~~~~~~~~~~~~~~~~~~~~~~~~~~~~~~~~~~~~~
    +2^{k-1}(2^{k}-1)\dfrac{2^{k}(2^{k}-1)^{3}}{2} \\
    &= \dfrac{(2^{k}+1)(2^{k}-1)(2^{k}-2)^{3}+2^{k-1}(2^{k}+1)(2^{k}-2)(2^{k}-3)^{3}+2^{2k-1}(2^{k}-1)^{4}}{2} \\
    &=\dfrac{2^{6k}-6\cdot2^{5k}+14\cdot2^{4k}-9\cdot2^{3k}-15\cdot2^{2k}+15\cdot2^{k}+8}{2}.
\end{align*}
We have 
\[
\dfrac{M_{1}(\mathcal{C}(G))}{|v(\mathcal{C}(G))|}=\dfrac{2^{5k}-4 \cdot 2^{4k} + 4 \cdot 2^{3k} +4 \cdot 2^{2k} -5 \cdot 2^{k} -4}{2^{3k}-2^{k}-1}
\] 
and
\[
\dfrac{M_{2}(\mathcal{C}(G))}{|e(\mathcal{C}(G))|}=\dfrac{2^{6k}-6\cdot2^{5k}+14\cdot2^{4k}-9\cdot2^{3k}-15\cdot2^{2k}+15\cdot2^{k}+8}{2^{4k}-2\cdot2^{3k}-2\cdot2^{2k}+3\cdot2^{k}+2}.
\]
Therefore,
\[
\dfrac{M_{2}(\mathcal{C}(G))}{|e(\mathcal{C}(G))|}- \dfrac{ M_{1}(\mathcal{C}(G))}{|v(\mathcal{C}(G))|}=\dfrac{2^{6k}(3\cdot2^{k}-11)+2^{k}(8\cdot2^{4k}-6\cdot2^{2k}-1)+2^{2k}(8\cdot2^{2k}-1)}{2^{5k}(2^{2k}-2\cdot2^{k}-3)+(4\cdot2^{4k}-2^{2k}-2)+2^{k}(6\cdot2^{2k}-5)} := \dfrac{f(k)}{g(k)}
\]
    For $ k \geq 2 $, we have $2\cdot2^{2k}(4\cdot2^{2k}-3)>1, 2^{k}(2^{k}-2)>3 $ and $2^{2k}(4\cdot2^{2k}-1)>2. $ Therefore, $\frac{f(k)}{g(k)}>0$ and so $\frac{M_{2}(\mathcal{C}(G))}{|e(\mathcal{C}(G))|} > \frac{M_{1}(\mathcal{C}(G))}{|v(\mathcal{C}(G))|}$.
    
    Using Theorem \ref{using commuting finding complement} we have
\begin{align*}
    M_{1}(\mathcal{NC}(G))&=(2^{3k}-2^{k}-1)(2^{3k}-2^{k}-2)^{2} \\
    &~~~~~~~~~~~~~-4 \cdot (2^{3k}-2^{k}-2) \dfrac{(2^{4k}-2\cdot2^{3k}-2\cdot2^{2k}+3\cdot2^{k}+2)}{2} \\
    &~~~~~~~~~~~~~~+(2^{5k}-4 \cdot 2^{4k} + 4 \cdot 2^{3k} +4 \cdot 2^{2k} -5 \cdot 2^{k} -4)  \\
    &= 2^{9k}-5\cdot2^{7k}-2^{6k}+9\cdot2^{5k}-5\cdot2^{3k}-3\cdot2^{2k}+3\cdot2^{k}
\end{align*}
and 
\begin{align*}
     M_{2}(\mathcal{NC}(G))&=\dfrac{(2^{3k}-2^{k}-1)(2^{3k}-2^{k}-2)^{3}}{2}+2\cdot \dfrac{(2^{4k}-2\cdot2^{3k}-2\cdot2^{2k}+3\cdot2^{k}+2)^{2}}{4} \\
     &~~~~~~~~~~~~~~~~~~~~~~~~~~~~~~~~~~~~~-3\cdot \dfrac{(2^{4k}-2\cdot2^{3k}-2\cdot2^{2k}+3\cdot2^{k}+2)}{2}(2^{3k}-2^{k}-2)^{2} \\
     &~~~~~~~~~~~~~~~~~~~~~~~~+(2^{3k}-2^{k}-1-\dfrac{3}{2})(2^{5k}-4 \cdot 2^{4k} + 4 \cdot 2^{3k} +4 \cdot 2^{2k} -5 \cdot 2^{k} -4) \\
    &~~~~~~~~~~~~~~~~~~~~~~~~~~~~~~~~- \dfrac{(2^{6k}-6\cdot2^{5k}+14\cdot2^{4k}-9\cdot2^{3k}-15\cdot2^{2k}+15\cdot2^{k}+8)}{2} \\
     &= \dfrac{1}{2}(2^{12k}-7\cdot2^{10k}-2^{9k}+21\cdot2^{8k}-26\cdot2^{6k}-2\cdot2^{5k}+15\cdot2^{4k}+3\cdot2^{3k} \\
     &~~~~~~~~~~~~~~~~~~~~~~~~~~~~~~~~~~~~~~~~~~~~~~~~~~~~~~~~~~~~~~~~~~~~~~~~~~~~~~~~~~~~~~~~~~~~~~+6\cdot2^{2k}-8\cdot2^{k}).
\end{align*}
    Also, $ |v(\mathcal{NC}(G))|=2^{3k}-2^{k}-1$ and $ |e( \mathcal{NC}(G))| = \binom{2^{3k}-2^{k}-1}{2}-|e(\mathcal{C}(G))|=\frac{1}{2}(2^{6k}-3\cdot2^{4k}-2^{3k}+3\cdot2^{2k}) $. We have
\[
\dfrac{M_{1}(\mathcal{NC}(G))}{|v( \mathcal{NC}(G) )|} = \dfrac{2^{9k}-5\cdot2^{7k}-2^{6k}+9\cdot2^{5k}-5\cdot2^{3k}-3\cdot2^{2k}+3\cdot2^{k}}{2^{3k}-2^{k}-1}
\]
and
\begin{align*}
    &\dfrac{M_{2}(\mathcal{NC}(G))}{|e(\mathcal{NC}(G))|} \\
    & ~~~~=  \dfrac{2^{12k}-7\cdot2^{10k}-2^{9k}+21\cdot2^{8k}-26\cdot2^{6k}-2\cdot2^{5k}+15\cdot2^{4k}+3\cdot2^{3k}+6\cdot2^{2k}-8\cdot2^{k}}{2^{6k}-3\cdot2^{4k}-2^{3k}+3\cdot2^{2k}}.    
\end{align*}
    As such 
\begin{align*}
    &\dfrac{M_{2}(\mathcal{NC}(G))}{|e( \mathcal{NC}(G) )|} - \dfrac{M_{1}(\mathcal{NC}(G))}{|v( \mathcal{NC}(G) )|} \\
    &~~~~= \dfrac{2^{7k}(2^{5k}-8\cdot2^{k}-4)+2^{4k}(17\cdot2^{2k}-14)+2^{3k}(14\cdot2^{2k}-18)+2\cdot2^{2k}+8\cdot2^{k}}{2^{6k}(2^{3k}-4\cdot2^{k}-2)+2\cdot2^{3k}(3\cdot2^{2k}-1)+2^{2k}(4\cdot2^{2k}-3)} := \dfrac{f(k)}{g(k)}.
\end{align*}
    For $ k \geq 2 $, we have $2^{k}(2^{4k}-8)>4 $ and $2^{k}(2^{2k}-4)>2. $ Therefore, $\frac{f(k)}{g(k)}>0$ and so $\frac{M_{2}(\mathcal{NC}(G))}{|e(\mathcal{NC}(G) )|} > \frac{M_{1}(\mathcal{NC}(G))}{|v( \mathcal{NC}(G) )|}.$
\end{proof}

We conclude this section with the following remark.

\begin{remark}\label{remark_for_groups}
The results of this section show that  Conjecture \ref{Conj} holds for commuting and non-commuting graphs of
\begin{enumerate}
\item the groups $D_{2m}, Q_{4n}, QD_{2^{n}}$, $V_{8n}$, $SD_{8n}$, $U_{6n}$, $M_{2mn}$, $S_{z}(2)$, $A(n,\nu)$, $A(n, p)$, $GL(2, q)$ and $PSL(2, 2^{k})$.
\item the non-abelian group of order $pq$,  where $p$ and $q$ are primes such that $p|q-1$.  
\item	the groups $G$  such that $\frac{G}{Z(G)} \cong D_{2m}$, $\mathbb{Z}_{p} \times \mathbb{Z}_{p}$ or $S_{z}(2)$. 
\end{enumerate} 
\end{remark}
\section{A few consequences}\label{consequences}
In this section we discuss the following consequences of the results obtained in Section 2.
\begin{theorem}
    Let $G$ be a finite non-abelian group and $|Z(G)|=n$.
    \begin{enumerate}
        \item If $G$ is 4-centralizer then
        $ M_{1}(\mathcal{C}(G))=3n(n-1)^{2}$, $M_{2}(\mathcal{C}(G))=\frac{3n(n-1)^{3}}{2}$, $M_{1}(\mathcal{NC}(G))=12n^{3}$ and $M_{2}(\mathcal{NC}(G))=18n^{4}$. \label{4-cent}
        \item If $G$ is 5-centralizer then $M_{1}(\mathcal{C}(G)) \in \{8n(2n-1)^{2}, 2n(2n-1)^{2}+3n(n-1)^{2}\}$, $M_{2}(\mathcal{C}(G)) \in \Big\{4n(2n-1)^{3}, \frac{1}{2}(2n(2n-1)^{2}+3n(n-1)^{2})\Big\}$, $M_{1}(\mathcal{NC}(G)) \in \{288n^{3}, 66n^{3}\}$ and $M_{2}(\mathcal{NC}(G)) \in \{1152n^{4}, 120n^{4}\}$. \label{5-cent}
        \item If $G$ is a $(p+2)$-centralizer $p$-group then
        $M_{1}(\mathcal{C}(G))=(pn-n)(p+1)(pn-n-1)^{2}$, $M_{2}(\mathcal{C}(G))=\frac{1}{2}(p+1)(pn-n)(pn-n-1)^{3}$, $M_{1}(\mathcal{NC}(G))$ $=(p+1)(pn-n)(p^{4}n^{2}-2p^{3}n^{2}+p^{2}n^{2}) $ and 
	$M_{2}(\mathcal{NC}(G))=\frac{1}{2}(p+1)^{2}(pn-n)^{2} 
        (p^{4}n^{2}-2p^{3}n^{2}+p^{2}n^{2})$.
        \item If \,$\{x_1, x_2, \ldots, x_r\}$ be a set of pairwise non-commuting elements of $G$ having maximal size, then for $r=3$, $M_{1}(\mathcal{C}(G))=3n(n-1)^{2}$, $M_{2}(\mathcal{C}(G))=\frac{3n(n-1)^{3}}{2}$, $M_{1}(\mathcal{NC}(G))=12n^{3}$ and $M_{2}(\mathcal{NC}(G))=18n^{4}$ and for $r=4$, $M_{1}(\mathcal{C}(G)) \in \{8n(2n-1)^{2}, 2n(2n-1)^{2}+3n(n-1)^{2}\}$, $M_{2}(\mathcal{C}(G)) \in \Big\{4n(2n-1)^{3}, \frac{1}{2}(2n(2n-1)^{2}+3n(n-1)^{2})\Big\}$, $M_{1}(\mathcal{NC}(G)) \in \{288n^{3}, 66n^{3}\}$ and $M_{2}(\mathcal{NC}(G)) \in \{1152n^{4}, 120n^{4}\}$.
    \end{enumerate}
    Further, $\frac{M_{2}(\Gamma(G))}{|e(\Gamma(G))|} \geq \frac{M_{1}(\Gamma(G))}{|v(\Gamma(G))|}$, where $\Gamma(G) = \mathcal{C}(G)$ or $\mathcal{NC}(G)$ in all the above cases.
\end{theorem}
\begin{proof}
    \begin{enumerate}
        \item By Theorem 2 of \cite{BS-MM-1994} we have that $\frac{G}{Z(G)} \cong \mathbb{Z}_{2} \times \mathbb{Z}_{2}$ when $G$ is 4-centralizer. Therefore, using Theorem \ref{thm-ZpXZp} and considering $p=2$ we get the required expressions for $M_{1}(\mathcal{C}(G)), M_{2}(\mathcal{C}(G)), M_{1}(\mathcal{NC}(G))$ and $M_{2}(\mathcal{NC}(G))$.
        \item By Theorem 4 of \cite{BS-MM-1994} we have that $\frac{G}{Z(G)} \cong \mathbb{Z}_{3} \times \mathbb{Z}_{3}$ or $D_6$ when $G$ is 5-centralizer. Therefore, using Theorem \ref{thm-ZpXZp} and Theorem \ref{thm-G/Z(G)=D_2m} and considering $p=3$ and $m=3$ respectively, we get the required expressions for $M_{1}(\mathcal{C}(G)), M_{2}(\mathcal{C}(G)), M_{1}(\mathcal{NC}(G))$ and $M_{2}(\mathcal{NC}(G))$.
        \item By Lemma 2.7 of \cite{AR-AC-2000} we have that $\frac{G}{Z(G)} \cong \mathbb{Z}_{p} \times \mathbb{Z}_{p}$ when $G$ is a $(p+2)$-centralizer $p$-group. Therefore, by Theorem \ref{thm-ZpXZp} we get the required expressions for $M_{1}(\mathcal{C}(G)), M_{2}(\mathcal{C}(G)),$ $M_{1}(\mathcal{NC}(G))$ and $M_{2}(\mathcal{NC}(G))$.
        \item By Lemma 2.4 of \cite{AJH-HJM-2007}, we have that $G$ is a $4$-centralizer or a $5$-centralizer group according as $r=3$ or $4$ if $\{x_1, x_2, \ldots, x_r\}$ is a set of pairwise non-commuting elements of $G$ having maximal size. Therefore, by parts \eqref{4-cent} and \eqref{5-cent} we get the desired expressions for $M_{1}(\mathcal{C}(G)), M_{2}(\mathcal{C}(G)), M_{1}(\mathcal{NC}(G))$ and $M_{2}(\mathcal{NC}(G))$.
    \end{enumerate}
    Also, by Theorem \ref{thm-ZpXZp} and Theorem \ref{thm-G/Z(G)=D_2m} we have $\frac{M_{2}(\Gamma(G))}{|e(\Gamma(G))|} \geq \frac{M_{1}(\Gamma(G))}{|v(\Gamma(G))|}$, where $\Gamma(G) = \mathcal{C}(G)$ or $\mathcal{NC}(G)$, in all the above cases.
\end{proof} 
\begin{theorem}
    Let $G$ be a finite non-abelian group with  $\Pr(G)$ as the commutativity degree of $G$ and $|Z(G)|=n$.
    \begin{enumerate}
        \item If $p$ is the smallest prime divisor of $|G|$ and  $\Pr(G)=\frac{p^{2}+p-1}{p^{3}}$ then $M_{1}(\mathcal{C}(G))=(pn-n)(p+1)(pn-n-1)^{2}$, $M_{2}(\mathcal{C}(G))=\frac{1}{2}(p+1)(pn-n)(pn-n-1)^{3}$, $M_{1}(\mathcal{NC}(G))$ $=(p+1)(pn-n)(p^{4}n^{2}-2p^{3}n^{2}+p^{2}n^{2}) $ and 
	$M_{2}(\mathcal{NC}(G))=\frac{1}{2}(p+1)^{2}(pn-n)^{2} 
        (p^{4}n^{2}-2p^{3}n^{2}+p^{2}n^{2})$.
        \item If \,$\Pr(G) \in \Big\{\frac{5}{14}, \frac{2}{5}, \frac{11}{27}, \frac{1}{2}, \frac{7}{16}, \frac{5}{8}\Big\}$ then $M_{1}(\mathcal{C}(G)) \in \{ 6n(6n-1)^{2}+7n(n-1)^{2}, 4n(4n-1)^{2}+5n(n-1)^{2}, 3n(3n-1)^{2}+4n(n-1)^{2}, 2n(2n-1)^{2}+3n(n-1)^{2}, 3n(n-1)^{2}, 8n(2n-1)^{2}\}$, $M_{2}(\mathcal{C}(G)) \in \Big\{\frac{1}{2}(6n(6n-1)^{3}+7n(n-1)^{3}), \frac{1}{2}(4n(4n-1)^{3}+5n(n-1)^{3}), \frac{1}{2}(3n(3n-1)^{3}+4n(n-1)^{3}), \frac{1}{2}(2n(2n-1)^{3}+3n(n-1)^{3}), \frac{3}{2}(n(n-1)^{3}), 4n(2n-1)^{3}\Big\}$, $M_{1}(\mathcal{NC}(G)) \in \{1302n^{3}, 420n^{3}, 192n^{3},66n^{3}, 12n^{3}, 288n^{3}\} $ and $M_{2}(\mathcal{NC}(G)) \in \{6552n^{4}, 1440n^{4}, 504n^{4},$ \\ $120n^{4}, 18n^{4}, 1152n^{4}\}$.
        \end{enumerate}
    Further, $\frac{M_{2}(\Gamma(G))}{|e(\Gamma(G))|} \geq \frac{M_{1}(\Gamma(G))}{|v(\Gamma(G))|}$, where $\Gamma(G) = \mathcal{C}(G)$ or $\mathcal{NC}(G)$ in both the above cases.
\end{theorem}
\begin{proof}
    \begin{enumerate}
        \item By Theorem 3 of \cite{DM-TMG-1974} we have that $\frac{G}{Z(G)} \cong \mathbb{Z}_{p} \times \mathbb{Z}_{p} $ if and only if $p$ is the smallest divisor of $|G|$ and  Pr$(G)=\frac{p^{2}+p-1}{p^{3}}$. Therefore, by Theorem \ref{thm-ZpXZp} we get the desired expressions for $M_{1}(\mathcal{C}(G)), M_{2}(\mathcal{C}(G)),$ $M_{1}(\mathcal{NC}(G))$ and $M_{2}(\mathcal{NC}(G))$.
        \item If $\Pr(G) \in \Big\{\frac{5}{14}, \frac{2}{5}, \frac{11}{27}, \frac{1}{2}, \frac{7}{16}, \frac{5}{8}\Big\}$ then by [\cite{DJR-PJM-1979}, pp. 246] and [\cite{RKN-AMS-2013}, pp. 451], we have $\frac{G}{Z(G)}$ is isomorphic to either $D_{14}, D_{10}, D_8, D_6, \mathbb{Z}_{2} \times \mathbb{Z}_{2}$ or $\mathbb{Z}_{3} \times \mathbb{Z}_{3}$. Therefore, by Theorem \ref{thm-G/Z(G)=D_2m} and Theorem \ref{thm-ZpXZp} we get the desired expressions for $M_{1}(\mathcal{C}(G)), M_{2}(\mathcal{C}(G)),$ $M_{1}(\mathcal{NC}(G))$ and $M_{2}(\mathcal{NC}(G))$.
    \end{enumerate}
    Also, by Theorem \ref{thm-G/Z(G)=D_2m} and Theorem \ref{thm-ZpXZp}, we have $\frac{M_{2}(\Gamma(G))}{|e(\Gamma(G))|} \geq \frac{M_{1}(\Gamma(G))}{|v(\Gamma(G))|}$, where $\Gamma(G) = \mathcal{C}(G)$ or $\mathcal{NC}(G)$, in both the above cases.
\end{proof}
\begin{theorem}\label{C(G)_is_planar}
    Let $G$ be a finite non-abelian group. If \, $\mathcal{C}(G)$ is planar, then $\frac{M_{2}(\Gamma(G))}{|e(\Gamma(G))|} \geq \frac{M_{1}(\Gamma(G))}{|v(\Gamma(G))|}$, where $\Gamma(G) = \mathcal{C}(G)$ or $\mathcal{NC}(G)$.
\end{theorem}
\begin{proof}
    By Theorem 2.2 of \cite{B3} we have that $\mathcal{C}(G)$ is planar if and only if $G$ is isomorphic to either $D_{6}, D_{8}, D_{10}, D_{12}, Q_{8}, Q_{12}, \mathbb{Z}_{2} \times D_{8}, \mathbb{Z}_{2} \times Q_{8}, \mathcal{M}_{16}, \mathbb{Z}_{4} \rtimes \mathbb{Z}_{4}, D_{8}*\mathbb{Z}_{4}, SG(16,3), A_{4}, A_{5}, S_{4}$, $SL(2,3) $ or $Sz(2).$ 	If $G \cong D_{6}, D_{8}, D_{10}, D_{12}, Q_{8},  Q_{12} \text{ or } Sz(2) $, then by Theorem \ref{thm-D2m}, Corollary \ref{cor-Q4n} and Corollary \ref{cor-Sz2} we have $\frac{M_{2}(\Gamma(G))}{|e(\Gamma(G))|} \geq \frac{M_{1}(\Gamma(G))}{|v(\Gamma(G))|}$, where $\Gamma(G) = \mathcal{C}(G)$ or $\mathcal{NC}(G)$.
	
	If $G \cong \mathbb{Z}_{2} \times D_{8}, \mathbb{Z}_{2} \times Q_{8}, \mathcal{M}_{16}, \mathbb{Z}_{4} \times \mathbb{Z}_{4}, D_{8}*\mathbb{Z}_{4} \text{ or } SG(16,3) $, then
	$\frac{G}{Z(G)} \cong \mathbb{Z}_{2} \times \mathbb{Z}_{2}$. Therefore, by Theorem \ref{thm-ZpXZp}, we have $\frac{M_{2}(\Gamma(G))}{|e(\Gamma(G))|} \geq \frac{M_{1}(\Gamma(G))}{|v(\Gamma(G))|}$, where $\Gamma(G) = \mathcal{C}(G)$ or $\mathcal{NC}(G)$.
	
    If $G \cong A_{4} $ then $\mathcal{C}(G)= K_{3} \sqcup 4K_{2}$. As such, $|v(\mathcal{C}(G))|=11, |e(\mathcal{C}(G))|= 7, M_{1}(\mathcal{C}(G))= 3(3-1)^{2} + 4 \cdot2 (2-1)^{2} = 20 \text{ and } M_{2}(\mathcal{C}(G))= 3 \cdot \frac{(3-1)^{3}}{2}+ 4 \cdot \frac{2(2-1)^{3}}{2} = 16. $ Therefore, 
\[
\frac{M_{2}(\mathcal{C}(G))}{|e(\mathcal{C}(G))|}=\dfrac{16}{7} > \dfrac{20}{11} = 
\frac{M_{1}(\mathcal{C}(G))}{|v(\mathcal{C}(G))|}.
\]
    Also, $|e(\mathcal{NC}(G))|= 48, M_{1}(\mathcal{NC}(G))= 11(11-1)^{2} - 4\cdot7(11-1)+20 = 840 \text{ and } M_{2}(\mathcal{NC}(G))= \frac{11(11-1)^{3}}{2}+2\cdot7^{2}-3\cdot7(11-1)^{2}+(11-\frac{3}{2})20-16 = 3672$. Therefore, 
\[
\frac{M_{2}(\mathcal{NC}(G))}{|e(\mathcal{NC}(G))|}=76.5 > \frac{840}{11} = 
\frac{M_{1}(\mathcal{NC}(G))}{|v(\mathcal{NC}(G))|}. 
\]
	
If $G \cong SL(2,3) $ then $\mathcal{C}(G)= 3K_{2} \sqcup 4K_{4}$.  As such, $|v(\mathcal{C}(G))|=22, |e(\mathcal{C}(G))|= 27, M_{1}(\mathcal{C}(G))= 3\cdot2(2-1)^{2} + 4 \cdot4 (4-1)^{2} = 150 \text{ and } M_{2}(\mathcal{C}(G))= 3 \cdot \frac{2(2-1)^{3}}{2}+ 4 \cdot \frac{4(4-1)^{3}}{2} = 219.$  Therefore,
\[
\frac{M_{2}(\mathcal{C}(G))}{|e(\mathcal{C}(G))|} = \dfrac{73}{9} > \dfrac{75}{11} =
\frac{M_{1}(\mathcal{C}(G))}{|v(\mathcal{C}(G))|}.
\]
    Also, $|e(\mathcal{NC}(G))|= 204, M_{1}(\mathcal{NC}(G))= 22(22-1)^{2} - 4\cdot27(22-1)+150 = 7584 \text{ and } M_{2}(\mathcal{NC}(G))= \frac{22(22-1)^{3}}{2}+2\cdot27^{2}-3\cdot27(22-1)^{2}+(22-\frac{3}{2})150-219 = 70464$. Therefore, 
\[
\frac{M_{2}(\mathcal{NC}(G))}{|e(\mathcal{NC}(G))|}=\frac{70464}{204} > \frac{7584}{22} = 
\frac{M_{1}(\mathcal{NC}(G))}{|v(\mathcal{NC}(G))|}. 
\]
    If $G \cong A_{5} $ then by Theorem \ref{thm-PSL} we have $ \frac{M_{2}(\mathcal{C}(G))}{|e(\mathcal{C}(G))|} \geq \frac{M_{1}(\mathcal{C}(G))}{|v( \mathcal{C}(G) )|}$ and $\frac{M_{2}(\mathcal{NC}(G))}{|e(\mathcal{NC}(G))|}\geq \frac{M_{1}(\mathcal{NC}(G))}{|v(\mathcal{NC}(G))|} $ since $A_{5} \cong PSL(2,4)$.
	
The commuting graph of $S_{4}$ is given by
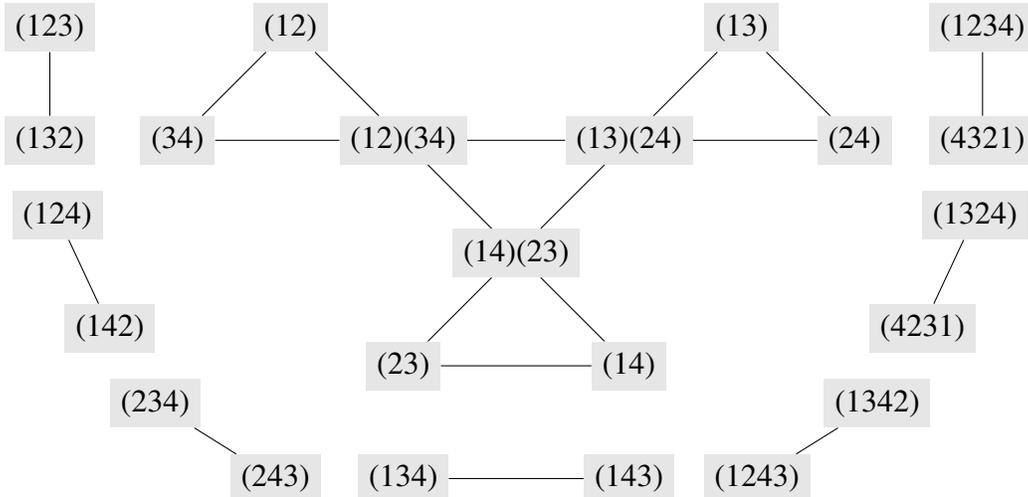
\begin{equation*}
	\begin{tikzpicture}
		\tikzstyle{vertex}=[fill=black!10]
		\node[vertex](a) at (0,0) {(12)(34)};
		\node[vertex](b) at (3,0) {(13)(24)};
		\node[vertex](c) at (1.5,-1.5) {(14)(23)};
		\node[vertex](d) at (-3,0) {(34)};
		\node[vertex](e) at (-1.5,1.5) {(12)};
		\node[vertex](f) at (4.5,1.5) {(13)};
		\node[vertex](g) at (6,0) {(24)};
		\node[vertex](h) at (0,-3) {(23)};
		\node[vertex](i) at (3,-3) {(14)};
		\node[vertex](j) at (-4.7,1.5) {(123)};
		\node[vertex](k) at (-4.7,0) {(132)};
		\node[vertex](l) at (-4.6,-1) {(124)};
		\node[vertex](m) at (-3.9,-2.5) {(142)};
		\node[vertex](n) at (-3.3,-3.5) {(234)};
		\node[vertex](o) at (-1.7,-4.5) {(243)};
		\node[vertex](p) at (0,-4.5) {(134)};
		\node[vertex](q) at (3,-4.5) {(143)};
		\node[vertex](r) at (7.7,1.5) {(1234)};
		\node[vertex](s) at (7.7,0) {(4321)};
		\node[vertex](t) at (7.6,-1) {(1324)};
		\node[vertex](u) at (6.9,-2.5) {(4231)};
		\node[vertex](v) at (6.3,-3.5) {(1342)};
		\node[vertex](w) at (4.7,-4.5) {(1243)};
		\path (a) edge (b)
		(b) edge (c)
		(a) edge (c)
		(a) edge (d)
		(a) edge (e)
		(e) edge (d)
		(b) edge (f)
		(g) edge (f)
		(h) edge (c)
		(i) edge (c)
		(b) edge (g)
		(h) edge (i)
		(j) edge (k)
		(l) edge (m)
		(n) edge (o)
		(p) edge (q)
		(r) edge (s)
		(t) edge (u)
		(v) edge (w);
	\end{tikzpicture}
\end{equation*}
    Therefore, if $G \cong S_{4}$ then $|v(\mathcal{C}(G))|=23, |e(\mathcal{C}(G))|= 19$,  $M_{1}(\mathcal{C}(G))= 86 \text{ and } M_{2}(\mathcal{C}(G))= 115.$  Hence, 
\[
\frac{M_{2}(\mathcal{C}(G))}{|e(\mathcal{C}(G))|}=\frac{115}{19} > \frac{86}{23} = \frac{M_{1}(\mathcal{C}(G))}{|v(\mathcal{C}(G))|}.
\]
     Also, $|e(\mathcal{NC}(G))|= 234, M_{1}(\mathcal{NC}(G))= 23(23-1)^{2} - 4\cdot19(23-1)+86 = 9456 \text{ and } M_{2}(\mathcal{NC}(G))= \frac{23(23-1)^{3}}{2}+2\cdot19^{2}-3\cdot19(23-1)^{2}+(23-\frac{3}{2})86-115 = 97320$. Therefore, 
\[
\frac{M_{2}(\mathcal{NC}(G))}{|e(\mathcal{NC}(G))|}=\frac{97320}{234} > \frac{9546}{23} = 
\frac{M_{1}(\mathcal{NC}(G))}{|v(\mathcal{NC}(G))|}. 
\]
  This completes the proof.
\end{proof}
\begin{theorem}\label{C(G)_is_toroidal}
Let $G$ be a finite non-abelian group. If $\mathcal{C}(G)$ is toroidal, then $\frac{M_{2}(\Gamma(G))}{|e(\Gamma(G))|} \geq \frac{M_{1}(\Gamma(G))}{|v(\Gamma(G))|}$, where $\Gamma(G) = \mathcal{C}(G)$ or $\mathcal{NC}(G)$.
\end{theorem} 
\begin{proof}
By Theorem 3.3 of \cite{B1} we have $\mathcal{C}(G)$ is toroidal if and only if $G$ is isomorphic to either $D_{14}, D_{16}, Q_{16}, QD_{16}, D_{6} \times \mathbb{Z}_{3}, A_{4} \times \mathbb{Z}_{2} \text{ or } \mathbb{Z}_{7} \rtimes \mathbb{Z}_{3}.$	If $G \cong D_{14}, D_{16}, Q_{16} \text{ or } QD_{16} $ then, by Theorem \ref{thm-D2m}, Corollary \ref{cor-Q4n} and Corollary \ref{cor-QD2n}, we have $\frac{M_{2}(\Gamma(G))}{|e(\Gamma(G))|} \geq \frac{M_{1}(\Gamma(G))}{|v(\Gamma(G))|}$, where $\Gamma(G) = \mathcal{C}(G)$ or $\mathcal{NC}(G)$.
If $G \cong \mathbb{Z}_{7} \rtimes \mathbb{Z}_{3}$ then $G$ is a group of order $pq$, where $p$ and $q$ are primes with $p|q-1$. Therefore, by Theorem \ref{thm-pq} we have $\frac{M_{2}(\Gamma(G))}{|e(\Gamma(G))|} \geq \frac{M_{1}(\Gamma(G))}{|v(\Gamma(G))|}$, where $\Gamma(G) = \mathcal{C}(G)$ or $\mathcal{NC}(G)$.
	
Note that	$D_{6} = \langle a, b: a^{3}=b^{2}=1, bab^{-1}=a^{-1} \rangle$ is an abelian centralizer group with center $Z(D_{6}) = \{1\}$ and  $C_{D_{6}}(a) = \lbrace 1, a, a^{2}\rbrace$, \, $C_{D_{6}}(ab) = \lbrace 1, ab \rbrace$,  $C_{D_{6}}(a^{2}b) = \lbrace 1, a^{2}b \rbrace$  and $C_{D_{6}}(b)= \lbrace 1, b \rbrace$ are the distinct centralizers of its non-central elements. Therefore, $D_{6} \times \mathbb{Z}_{3}$ is also an abelian centralizer group with center  $Z(D_{6} \times \mathbb{Z}_{3}) = \{1, a^2\} \times \mathbb{Z}_{3}$  and  $\lbrace 1, a, a^{2} \rbrace \times \mathbb{Z}_{3}$, \, $\lbrace 1, ab \rbrace \times \mathbb{Z}_{3}$,  $\lbrace 1, a^{2}b \rbrace \times \mathbb{Z}_{3}$  and $\lbrace 1, b \rbrace \times \mathbb{Z}_{3}$ are the distinct centralizers of non-central elements of $D_{6} \times \mathbb{Z}_{3}$.
Hence, if $G \cong D_{6} \times \mathbb{Z}_{3}$ then,  by Lemma 2.1 of \cite{B1}, we have $\mathcal{C}(G)= K_{6} \sqcup 3K_{3}$. As such,  $|v(\mathcal{C}(G))| = 15, |e(\mathcal{C}(G))| = 24, M_{1}(\mathcal{C}(G))= 6\cdot(6-1)^{2} + 3 \cdot3 (3-1)^{2} = 186 \text{ and } M_{2}(\mathcal{C}(G))= 6 \cdot \frac{(6-1)^{3}}{2}+ 3 \cdot \frac{3(3-1)^{3}}{2} = 411.$  Therefore, 
\[
\frac{M_{2}(\mathcal{C}(G))}{|e(\mathcal{C}(G))|} = 17.125 > 12.4 = \frac{M_{1}(\mathcal{C}(G))}{|v(\mathcal{C}(G))|}.
\]
   Also, \, $|e(\mathcal{NC}(G))|= 81, M_{1}(\mathcal{NC}(G))= 15(15-1)^{2} - 4\cdot24(15-1)+186 = 1782 \text{ and } M_{2}(\mathcal{NC}(G))= \frac{15(15-1)^{3}}{2}+2\cdot24^{2}-3\cdot24(15-1)^{2}+(15-\frac{3}{2})186-411 = 9720$. Therefore, 
\[
\frac{M_{2}(\mathcal{NC}(G))}{|e(\mathcal{NC}(G))|}=120 > 118.8 = 
\frac{M_{1}(\mathcal{NC}(G))}{|v(\mathcal{NC}(G))|}. 
\]  

    We have	$A_{4} = \langle a, b: a^{2}=b^{3}=(ab)^{3} = 1 \rangle$ is an abelian centralizer group with center $Z(A_4) = \{1\}$ and $C_{A_{4}}(a)=\lbrace 1, a, bab^{2}, b^{2}ab \rbrace$, \, $C_{A_{4}}(ab)=\lbrace 1, ab, b^{2}a \rbrace$, \, $C_{A_{4}}(aba)= \lbrace 1, aba, bab \rbrace$, \, $C_{A_{4}}(b)= \lbrace 1, b, b^{2} \rbrace$ and $C_{A_{4}}(ba)= \lbrace 1, ba, ab^{2} \rbrace$ are the distinct centralizers of its non-central elements. Therefore, $A_{4} \times \mathbb{Z}_{2}$ is also an abelian centralizer group with center $Z(A_{4} \times \mathbb{Z}_{2}) = \{1\} \times \mathbb{Z}_{2}$ and $\lbrace 1, a, bab^{2}, b^{2}ab \rbrace \times \mathbb{Z}_{2}$, \, $\lbrace 1, ab, b^{2}a \rbrace \times \mathbb{Z}_{2}$, \, $\lbrace 1, aba, bab \rbrace \times \mathbb{Z}_{2}$, \, $\lbrace 1, b, b^{2} \rbrace \times \mathbb{Z}_{2}$ and $\lbrace 1, ba, ab^{2} \rbrace \times \mathbb{Z}_{2}$ are the distinct centralizers of  non-central elements of $A_{4} \times \mathbb{Z}_{2}$. Hence, if $G \cong A_{4} \times \mathbb{Z}_{2}$ then, by Lemma 2.1 of \cite{B1}, we have $\mathcal{C}(G)= K_{6} \sqcup 4K_{4}$. As such,  $|v(\mathcal{C}(G))|=22, |e(\mathcal{C}(G))|= 39, M_{1}(\mathcal{C}(G))= 6\cdot(6-1)^{2} + 4 \cdot4 (4-1)^{2} = 294 \text{ and } M_{2}(\mathcal{C}(G))= 6 \cdot \frac{(6-1)^{3}}{2}+ 4 \cdot \frac{4(4-1)^{3}}{2} = 591.$  Therefore, 
\[
\frac{M_{2}(\mathcal{C}(G))}{|e(\mathcal{C}(G))|}=\frac{197}{13} > \frac{147}{11} = \frac{M_{1}(\mathcal{C}(G))}{|v(\mathcal{C}(G))|}.
\]
    Also, \, $|e(\mathcal{NC}(G))|= 192, M_{1}(\mathcal{NC}(G))= 22(22-1)^{2} - 4\cdot39(22-1)+294 = 6720 \text{ and } M_{2}(\mathcal{NC}(G))= \frac{22(22-1)^{3}}{2}+2\cdot39^{2}-3\cdot39(22-1)^{2}+(22-\frac{3}{2})294-591 = 58752$. Therefore, 
\[
\frac{M_{2}(\mathcal{NC}(G))}{|e(\mathcal{NC}(G))|}=306 > \frac{3360}{11}= 
\frac{M_{1}(\mathcal{NC}(G))}{|v(\mathcal{NC}(G))|}. 
\]
    This completes the proof.
\end{proof}
\begin{theorem}\label{NC(G)_is_planar}
    Let $G$ be a finite non-abelian group. If $\mathcal{NC}(G)$ is planar, then $\frac{M_{2}(\Gamma(G))}{|e(\Gamma(G))|} \geq \frac{M_{1}(\Gamma(G))}{|v(\Gamma(G))|}$, where $\Gamma(G) = \mathcal{C}(G)$ or $\mathcal{NC}(G)$. 
\end{theorem}
\begin{proof}
    If $\mathcal{NC}(G)$ is planar then by Proposition 2.3 of \cite{B5} we have that $ G $ is isomorphic to either $D_{6}, D_{8}$ or $Q_{8}$. In any of the above mentioned cases, we get $\frac{M_{2}(\Gamma(G))}{|e(\Gamma(G))|} \geq \frac{M_{1}(\Gamma(G))}{|v(\Gamma(G))|}$, where $\Gamma(G) = \mathcal{C}(G)$ or $\mathcal{NC}(G)$ by Theorem \ref{thm-D2m} and Corollary \ref{cor-Q4n}. 
\end{proof}
We conclude this section with the following corollary.
\begin{cor}
    Let $G$ be a finite non-abelian group.
\begin{enumerate}
\item If \, $\mathcal{C}(G)$ is planar then $M_{1}(\mathcal{C}(G)) \in \{2, 6, 20, 36, 42, 86, 96, 108, 150,  296\}$, $M_{2}(\mathcal{C}(G)) \in \{1, 3, 16, 54, 57, 114, 115, 162, 219, 394\}$, $M_{1}(\mathcal{NC}(G)) \in \{66, 96, 420, 528, 768, 840, 4740,$ $7584, 9546, 184988\} $ and $M_{2}(\mathcal{NC}(G)) \in \{120, 192, 1440, 1920, 3672, 4608, 37440, 70464,$ $97320, 5223424\}$.
\item If \, $\mathcal{C}(G)$ is toroidal then $M_{1}(\mathcal{C}(G)) \in \{150, 158, 164, 186, 294\}$, $M_{2}(\mathcal{C}(G)) \in \{375, 379, 382,$ $411, 591\}$, $M_{1}(\mathcal{NC}(G)) \in \{1302, 1536, 1782, 6299, 6720\} $ and $M_{2}(\mathcal{NC}(G)) \in \{6552, 8064,$ $9720, 58752, 76127\}$.
\item If $\mathcal{NC}(G)$ is planar then $M_{1}(\mathcal{C}(G)) \in \{2, 6\}$, $M_{2}(\mathcal{C}(G)) \in \{1, 3\}$, $M_{1}(\mathcal{NC}(G)) \in \{66, 96\} $ and $M_{2}(\mathcal{NC}(G)) \in \{120, 192\}$.
\end{enumerate} 
\end{cor}
\section{Concluding remark}
    As mentioned in Remark \ref{remark_for_groups}, we have found that the Conjecture \ref{Conj} holds for the commuting and non-commuting graphs of several families of finite groups. In Section \ref{consequences}, we have found that when a finite group satisfies certain conditions, its commuting and non-commuting graphs also satisfy Conjecture \ref{Conj}. 

     Also, using the following  GAP program, we have found that the commuting and non-commuting graphs of  finite non-abelian groups up to order $1000$
     satisfy Conjecture \ref{Conj}.
 \begin{verbatim}
LoadPackage("grape");
ComGraph:=function(G) 
 local vert,rel;
     if IsAbelian(G) then Error("Group must be non-abelian"); fi;
     vert:=Difference(G,Center(G));
     rel:={x,y}->x<>y and x*y=y*x;
     return Graph(Group(()),vert,{x,g}->x,rel,true);
 end;
 
HVCon:=function(Gr) 
local M1,M2,Grc;
   M1:=Sum(Vertices(Gr),v->VertexDegree(Gr,v)^2)/Size(Vertices(Gr));
   M2:=Sum(UndirectedEdges(Gr),
       e->VertexDegree(Gr,e[1])*VertexDegree(Gr,e[2]))/
                                             Size(UndirectedEdges(Gr));
   if M2<M1 then return false; fi;
   Grc:=ComplementGraph(Gr);
   M1:=Sum(Vertices(Grc),v->VertexDegree(Grc,v)^2)/Size(Vertices(Grc));
   M2:=Sum(UndirectedEdges(Grc),
       e->VertexDegree(Grc,e[1])*VertexDegree(Grc,e[2]))/
                                            Size(UndirectedEdges(Grc));
   if M2<M1 then return false; else return true; fi;
 end;
 for d in [1..1000] do
     Print(d,"\n");
     for id in [1..NrSmallGroups(d)] do
         G:=SmallGroup(d,id);
         if not IsAbelian(G) and not HVCon(ComGraph(G)) 
         then Print("found",[d,id],"\n"); fi;
     od;
 od;
\end{verbatim}

In view of above discussion,  we conclude this paper with the following conjecture.
      \begin{conj}
          Let $G$ be a finite non-abelian group. If $\Gamma(G)$ denotes the
           commuting or non-commuting graph of $G$, then $$\frac{M_{2}(\Gamma(G))}{|e(\Gamma(G))|} \geq \frac{M_{1}(\Gamma(G))}{|v(\Gamma(G))|}.$$
      \end{conj}

\section*{Acknowledgements}
The first author is thankful to Council of Scientific and Industrial Research  for the fellowship (File No. 09/0796(16521)/2023-EMR-I). The authors would like to thank Benjamin Sambale, Institut f\"ur Algebra, Zahlentheorie und Diskrete Mathematik, Leibniz Universit\"at Hannover, 30167 Hannover, Germany for helping with the GAP code.

\end{document}